\documentclass{article}
\usepackage{amsmath, amsfonts, amssymb, amsthm, mathrsfs}
\usepackage{tikz-cd}
\usepackage{xcolor}
\usepackage{ytableau}
\usepackage[backrefs,msc-links]{amsrefs}
\usepackage{hyperref} 
\usepackage[all]{xy}
\usepackage{caption}
\hypersetup{hidelinks}

\newtheorem{theorem}{Theorem}[section]
\newtheorem{proposition}[theorem]{Proposition}
\newtheorem{corollary}[theorem]{Corollary}
\newtheorem{lemma}[theorem]{Lemma}
\newtheorem{remark}[theorem]{Remark}
\newtheorem{definition}[theorem]{Definition}

\newtheorem{example}[theorem]{Example}

\def\loccitt{\emph{loc.\ cit.}}
\def\loccit{\emph{loc.\ cit.\ }}

\def\fsl{{\mathfrak{sl}}}
\def\fgl{{\mathfrak{gl}}}

\def\fZ{{\mathfrak{Z}}}

\def\BA{{\mathbb{A}}}
\def\BC{{\mathbb{C}}}

\def\BF{{\mathbb{F}}}
\def\BP{{\mathbb{P}}}
\def\BR{{\mathbb{R}}}
\def\BQ{{\mathbb{Q}}}
\def\BZ{{\mathbb{Z}}}

\def\CA{{\mathcal{A}}}
\def\CB{{\mathcal{B}}}

\def\CF{{\mathcal{F}}}
\def\CG{{\mathcal{G}}}

\def\CL{{\mathcal{L}}}
\def\CM{{\mathcal{M}}}

\def\CO{{\mathcal{O}}}

\def\CV{{\mathcal{V}}}

\def\sT{{\mathsf{T}}}
\def\sA{{\mathsf{A}}}

\def\rV{{\mathscr{V}}}

\def\ofZ{\overline{\fZ}}

\def\Hom{\textrm{Hom}}

\def\col{\textrm{col }}

\def\and{\textrm{ }\&\textrm{ }}

\def\op{\overline{p}}

\def\ofZ{\overline{\fZ}}

\def\glnhat{{U_q(\widehat{\fgl}_n)}}
\def\glnhatleq{{U_q^\leq(\widehat{\fgl}_n)}}
\def\glnhatgeq{{U_q^\geq(\widehat{\fgl}_n)}}
\def\uup{{U_q^+(\widehat{\fgl}_n)}}
\def\uum{{U_q^-(\widehat{\fgl}_n)}}
\def\slnhat{{U_q(\widehat{\fsl}_n)}}

\def\bd{{\mathbf{d}}}

\def\bu{{\mathbf{u}}}

\def\la{{\lambda}}
\def\sq{{\square}}
\def\bsq{{\blacksquare}}

\def\bla{{\boldsymbol{\la}}}
\def\bmu{{\boldsymbol{\mu}}}
\def\bnu{{\boldsymbol{\nu}}}

\def\bA{{\mathbf{A}}}
\def\bB{{\mathbf{B}}}

\def\bV{{\mathbf{V}}}
\def\bX{{\mathbf{X}}}
\def\bY{{\mathbf{Y}}}

\def\const{\mathrm{const}}

\newcommand{\hide}[1]{}
\newcommand{\todo}[1]{}

\title{Affine Laumon space and contragredient dual Verma module of $\glnhat$}
\author{Che Shen}
\date{}

\begin{document}

\maketitle

\begin{abstract}
We study the action of the quantum group $\glnhat$ on the equivariant K-theory of affine Laumon spaces. 
We show that, at any lowest weight away from the critical level, this can be identified with the contragredient dual Verma module of $\glnhat$, 
improving earlier results of Braverman-Finkelberg and Negu{\c{t}}. 
The proof uses a variant of stable envelopes first introduced by Maulik-Okounkov 
in the study of Nakajima quiver varieties.


\end{abstract}

\tableofcontents

\section{Introduction}
\subsection{Laumon space} \label{finite_case}
Laumon spaces parametrize flags of locally-free sheaves on $\mathbb{P}^1$:
\begin{equation} 
    \label{flagofsheaves}
    \CF_1 \subset \dots \subset \CF _{n-1} \subset \CO _{\BP^1}^{\oplus n}
\end{equation}
whose fibers near $\infty \in \BP^1$ match a fixed full flag of subspaces of $\BC^n$.
Laumon spaces are disconnected, with connected components indexed by vectors $\bd = (d_1,\dots,d_{n-1}) \in \BZ _{\geq 0} ^{n-1}$
that keep track of the degree of the sheaves in (\ref{flagofsheaves}).  
We denote these connected components by $\CM_\bd^\text{fin}$.\footnote{The superscript ``fin'' stands for finite, which reflects its relation to the finite-dimensional Lie algebra, whereas affine Laumon space relates to the affine Lie algebra.}
The torus $\sT = (\mathbb{C}^*)^{n+1}$ acts on $\CM_\bd^\text{fin} $ by scaling the fibers of $\CO_{\BP^1}^{\oplus n}$ as well as scaling the base $\mathbb{P}^1$.

The relation between Laumon spaces and representation theory of $\mathfrak{sl}_n$, 
as well as their affine analogs, was extensively studied in the literature, 
see e.g. \cites{feigin2011gelfand,nakajima2011handsaw,feigin2011yangians,tsymbaliuk2010quantum, neguct2018affine}. 
(See also \cites{bullimore2016vortices,hilburn2023bfn} in the context of supersymmetric gauge theories.)
In particular, it was proved in \cite{braverman2005finite} that the localized equivariant K-theory of Laumon space 
\[ 
    K ^\text{loc, fin}:= \bigoplus_{\bd \in \BZ _{\geq 0}^{n-1}}K_\sT(\CM_\bd^\text{fin} ) \otimes _{K_\sT(\text{pt})}\text{Frac}\, K_\sT(\text{pt})
\]
has a natural
action 
\footnote{\label{note_larger} In fact, a larger algebra action was constructed, see Remark \ref{larger_alg}.} 
of the quantum group $U_q(\fsl_n)$ and can be identified with the universal Verma module of $U_q(\fsl_n)$ over $\text{Frac}\, K_\sT(\text{pt})$, 
whose lowest weight corresponds to equivariant parameters $(u_1,\dots,u_n)$. It is natural to ask what module we get when specializing the parameters to a certain lowest weight 
\begin{equation}\label{specialize_hwt}
    u_i = q ^{a_i}, \; a_i \in \BC, \; i=1,\dots,n
\end{equation}
The answer is less obvious than it seems because certain coefficients might become 0 or $\infty$, making it no longer a Verma module. In fact, as a corollary of Theorem \ref{maintheorem}, we prove that the (non-localized) K-theory
$$K ^\text{int, fin}:= \bigoplus_{\bd \in \BZ _{\geq 0}^{n-1}}K_\sT(\CM_\bd^\text{fin} )$$ 
becomes the \textit{contragredient 
dual} of the Verma module of $U_q(\fsl_n)$ under the specialization (\ref{specialize_hwt}). 
A similar result was obtained in \cite{feigin2011gelfand} when considering (non-quantized) $U(\mathfrak{gl}_n)$ acting on the cohomology of Laumon space.

\subsection{Affine Laumon space} \label{aff_case}
There is an affine analogue of the above result. Consider the affine Laumon space $\CM_\bd $, 
see section \ref{aff_laumon} for the definition. \cites{tsymbaliuk2010quantum,neguct2018affine} constructed
 a geometric action of $\glnhat$ on the localized equivariant K-theory
\[ 
    K ^\text{loc} := \bigoplus_{\bd \in \BZ _{\geq 0}^{n}}K_\sT(\CM_\bd) \otimes _{K_\sT(\text{pt})}\text{Frac}\, K_\sT(\text{pt})
\]
and \cite{neguct2018affine} further proved that this makes $K ^\text{loc}$ the universal Verma module of $\glnhat$ over $\text{Frac}\, K_\sT(\text{pt})$. Similar to the finite case,  
we prove in corollary \ref{dual_verma_cor_aff} that 
\[ 
    K ^\text{int}:= \bigoplus_{\bd \in \BZ _{\geq 0}^{n}}K_\sT(\CM_\bd)
\]
becomes the contragredient dual Verma module of $\glnhat$ 
when specializing the equivariant parameters to a certain lowest weight. As we will explain in section \ref{action_almost_preserve},  the lowest weight must be away from the critical level due to a mild denominator in the $\glnhat$ action.

\subsection{Strategy of proof}
When $d_n=0$,  the affine Laumon space $\CM_\bd $ becomes the (finite) Laumon space $\CM^\text{fin} _{(d_1,\dots,d _{n-1})}$, 
so the results discussed in Section \ref{finite_case} could be a corollary of Section \ref{aff_case}. 
However, in this paper, we first present a proof for the finite Laumon space, which is simpler than the proof in the general case. The proof makes crucial use of the upper triangular property of the PBW basis, see Section \ref{proof_finite_case}. Some key ideas are reused in the proof for the affine case. 
 

For affine Laumon space, the upper triangular property of the PBW basis no longer holds and new methods are required. 
It turns out that a variant of stable envelope, first defined in \cite{maulik2012quantum} and further developed in \cites{okounkov2015lectures,okounkov2021inductive}, provides a way to ``restore'' the upper triangular property and plays a key role in the proof. In Section \ref{proof_aff_case}, we explain in detail how to use the properties of stable envelope and the rigidity argument in \cite{okounkov2015lectures} to compute their pairing with the PBW basis. 

In the affine case, there is an additional challenge that the correspondences used to define lowering operators are not proper, and thus introducing denominators in the action. To specialize equivariant parameters, we need to get good control of the denominators. This is done using a combination of algebraic and geometric arguments. It turns out that these denominators are very mild, see Proposition \ref{coef_almost_int}.

\subsection*{Acknowledgement}
I am grateful to my advisor Andrei Okounkov for many valuable discussions and suggestions. I would also like to thank Shaoyun Bai, Sam DeHority, Hunter Dinkins, Yixuan Li, Andrei Negu{\c{t}}, Leonid Rybnikov, Spencer Tamagni, Tianqing Zhu, and the anonymous referee for many useful comments and suggestions.

\section{Laumon spaces and quantum affine algebra}
This section collects some basic properties of affine Laumon space as well as the 
action of quantum affine algebra $\glnhat$ on the localized K-theory of it. 
We will closely follow the notations and figures in \cite{neguct2018affine}. 

\subsection{Affine Laumon space} \label{aff_laumon}
\subsubsection{Definition as a quiver variety}
Given an integer $n \geq 2$ and a degree $\bd \in \BZ_{\geq 0}^{n}$, 
affine Laumon spaces $\CM_\bd$ can be described in two equivalent ways: either as the moduli space of framed parabolic sheaves on $\BP^1 \times \BP^1$, or as a quiver variety, 
cf. \cite{neguct2018affine}, section 3. We will use the latter description in this paper.

Consider the following quiver:

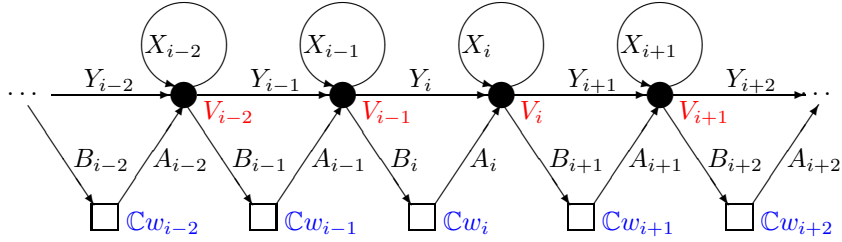
\begin{picture}(200,110)(40,-30)

    \put(43,31){\dots}
    \put(343,31){\dots}
    
    \put(60,31){\vector(1,0){45}}
    \put(72,34){$Y_{i-2}$}
    
    \put(115,31){\vector(1,0){50}}
    \put(135,34){$Y_{i-1}$}
    
    \put(175,31){\vector(1,0){50}}
    \put(195,34){$Y_{i}$}
    
    \put(235,31){\vector(1,0){50}}
    \put(255,34){$Y_{i+1}$}
    
    \put(295,31){\vector(1,0){50}}
    \put(315,34){$Y_{i+2}$}
    
    \put(110,31){\circle*{10}}
    \put(110,50){\circle{30}}
    \put(95,47){$X_{i-2}$}
    \put(117,22){$\color{red}{V_{i-2}}$}
    \put(102,36){\vector(4,-1){5}}
    
    \put(170,31){\circle*{10}}
    \put(170,50){\circle{30}}
    \put(155,47){$X_{i-1}$}
    \put(177,22){$\color{red}{V_{i-1}}$}
    \put(162,36){\vector(4,-1){5}}
    
    \put(230,31){\circle*{10}}
    \put(230,50){\circle{30}}
    \put(215,47){$X_{i}$}
    \put(237,22){$\color{red}{V_{i}}$}
    \put(222,36){\vector(4,-1){5}}
    
    \put(290,31){\circle*{10}}
    \put(290,50){\circle{30}}
    \put(275,47){$X_{i+1}$}
    \put(297,22){$\color{red}{V_{i+1}}$}
    \put(282,36){\vector(4,-1){5}}
    
    \put(75,-20){\line(1,0){10}}
    \put(75,-10){\line(1,0){10}}
    \put(75,-20){\line(0,1){10}}
    \put(85,-20){\line(0,1){10}}
    \put(88,-20){$\color{blue}{\BC w_{i-2}}$}
    
    \put(51,27){\vector(2,-3){25}}
    \put(85,-10){\vector(2,3){25}}
    \put(68,4){$B_{i-2}$}
    \put(98,4){$A_{i-2}$}
    
    \put(135,-20){\line(1,0){10}}
    \put(135,-10){\line(1,0){10}}
    \put(135,-20){\line(0,1){10}}
    \put(145,-20){\line(0,1){10}}
    \put(148,-20){$\color{blue}{\BC w_{i-1}}$}
    
    \put(111,27){\vector(2,-3){25}}
    \put(145,-10){\vector(2,3){25}}
    \put(128,4){$B_{i-1}$}
    \put(158,4){$A_{i-1}$}
    
    \put(195,-20){\line(1,0){10}}
    \put(195,-10){\line(1,0){10}}
    \put(195,-20){\line(0,1){10}}
    \put(205,-20){\line(0,1){10}}
    \put(208,-20){$\color{blue}{\BC w_{i}}$}
    
    \put(171,27){\vector(2,-3){25}}
    \put(205,-10){\vector(2,3){25}}
    \put(188,4){$B_{i}$}
    \put(218,4){$A_{i}$}
    
    \put(255,-20){\line(1,0){10}}
    \put(255,-10){\line(1,0){10}}
    \put(255,-20){\line(0,1){10}}
    \put(265,-20){\line(0,1){10}}
    \put(268,-20){$\color{blue}{\BC w_{i+1}}$}
    
    \put(231,27){\vector(2,-3){25}}
    \put(265,-10){\vector(2,3){25}}
    \put(248,4){$B_{i+1}$}
    \put(278,4){$A_{i+1}$}
    
    \put(315,-20){\line(1,0){10}}
    \put(315,-10){\line(1,0){10}}
    \put(315,-20){\line(0,1){10}}
    \put(325,-20){\line(0,1){10}}
    \put(328,-20){$\color{blue}{\BC w_{i+2}}$}
    
    \put(291,27){\vector(2,-3){25}}
    \put(325,-10){\vector(2,3){25}}
    \put(308,4){$B_{i+2}$}
    \put(338,4){$A_{i+2}$}
    
    \label{chainsaw_quiver}
\end{picture} 
%

\noindent Note that the quiver is meant to be cyclic, with a total of $n$ circle nodes $V_1,...,V_n$ and $n$ framing nodes $\mathbb{C}w_1,...,\mathbb{C}w_n$. All indices are understood mod $n$. For examples, the arrow $Y_1$ goes from $V_n$ to $V_1$; the arrow $B_1$ goes from $V_n$ to $\mathbb{C}w_1$.
Fix vector spaces $V_1, \dots, V_n$ of dimension $d_1, \dots, d_n$, consider the vector space of linear maps 
\begin{equation}
    \label{prequotient}
    M_\bd = \bigoplus_{i=1}^{n} \left(\Hom(V_i, V_i) \oplus \Hom(V_{i-1}, V_i) \oplus \Hom(W_i, V_i) \oplus \Hom(V_{i-1}, W_i) \right)\qquad
\end{equation}
where we set $V_{n+1} = V_1$ and $V_0 = V_n$. 
Elements of the vector space $M_\bd$ will be quadruples $(X_i, Y_i, A_i, B_i)_{1 \leq i \leq n}$. Consider the quadratic map 
\[ 
    M_\bd \stackrel{\nu}{\longrightarrow} \bigoplus_{i=1}^{n} \Hom(V_{i-1}, V_i)
\]
\[
\nu(X_i, Y_i, A_i, B_i)_{1 \leq i \leq n} = \bigoplus_{i=1}^{n} \left( X_i Y_i - Y_i X_{i-1} + A_i B_i \right)
\]
and define
\[ 
    \CM_\bd = \nu^{-1}(0)^{\text{stable}} / GL_\bd 
\]
where a point in $\nu^{-1}(0)$ is stable if the $V_i$'s are generated by $X$ and $Y$ acting on the images of $A$ maps.

There is a torus action by 
\begin{equation}
    \label{torus}
    \sT = (\mathbb{C}^*)^{n+2} = T_W \times \mathbb{C}_q^* \times \BC_p^*
\end{equation} 
Explicitly, if $(u_1, \dots, u_n, q, p)$ denote\footnote{$p$ was denoted $\overline{q}$ in \cite{neguct2018affine}} the characters of $\sT$, then $\sT$
scales the quadruple $(X_i, Y_i, A_i, B_i)$ with weight 
\[ 
    (q^{-2}, p^{-2 \delta_i^1}, u_i^{-2}, u_i^2 q^{-2} p^{-2 \delta_i^1})
\]
($p$ only acts on $Y_1$ and $B_1$.) 
Note that we are using a $2^{n+2}$-fold cover of the usual torus action to avoid square roots in formulas. \hide{I think Negut paper used the inverse of this action, which doesn't match with the formulas... even for $d=(1,0)$.}

Define
\[ 
    K^\text{int} := \bigoplus_{\bd \in \BZ_{\geq 0}^n} K_\bd^\text{int} := \bigoplus_{\bd \in \BZ_{\geq 0}^n} K_\sT(\CM_\bd)
\]
which is a module over
\begin{equation} \label{Zu}
    \BZ_\bu := \BZ[u_1^\pm, \dots, u_n^\pm, q^\pm, p^\pm] = K_\sT(\text{pt})
\end{equation}
We also define localized equivariant K-theory
\[ 
    K^\text{loc} := K^\text{int} \otimes_{\BZ_\bu} \BF_\bu
\]
where
\begin{equation}\label{Fu}
    \BF_\bu := \text{Frac}\, \BZ_\bu = \BQ(u_1, \dots, u_n, q, p)
\end{equation}

\subsubsection{Fixed points} \label{fixedpt} 
The fixed points are parametrized by $n$-tuples of 2d partitions $\bla = (\lambda^1, \dots, \lambda^n)$ 
where each $\lambda^i$ itself is a 2d partition. We will simply call $\bla$ a partition for brevity. (It should be clear from the context whether we are talking about
a single 2d partition or an $n$-tuple of 2d partitions.) For a box $\sq$ in the $(x, y)$ position of the partition $\lambda^k$, we say that $\sq$ has color $y + k \mod n$. (Note that $x, y$ start from 0.) Write $|\bla| = \bd$ where $\bd = (d_1, \dots, d_n) \in \BZ_{\geq 0}^n$ 
if $\bla$ has $d_i$ boxes of color $i$, $i = 1, 2, \dots, n$.

Fixed points of $\CM_\bd$ are parametrized by partitions $\bla$ with $|\bla| = \bd$. We will abuse notation to use $\bla$ to denote either a fixed point or a partition. The quadruple $(X,Y,A,B)$ corresponding to the fixed point $\bla$ has each vector space $V_i$ spanned by color $i$ boxes in $\bla$, i.e.
\[ 
    V_i = \bigoplus_{\sq \in \bla}^{\col \sq = i}\mathbb{C}\cdot v_\sq
\]
and for any $v_\sq \in V_i$,
\[ 
    X_i \cdot v_\sq = v _{\text{box directly to the right of } \sq}
\]
\[ 
    Y_{i+1} \cdot v_\sq = v _{\text{box directly above } \sq}
\]
\[ 
    A \cdot w_i = v _{\text{corner of }\lambda ^i}
\]

We will use the normalized fixed point basis
\[ 
    |\bla \rangle := \frac{[\CO_\bla]}{\Lambda^\bullet(T ^\vee_\bla  \CM_\bd)}
\]
where $[\CO_\bla]$ is the structure sheaf supported at the point $\bla$. With this notation, the equivariant localization formula becomes
\[ 
    c = \sum_{|\bla| = \bd} c|_\bla \cdot |\bla \rangle 
\]
where $c \in K_\sT(\CM_\bd)$.

\subsubsection{Local coordinates} \label{local_coord}
For any fixed point $\bla$, there is an open subset $U_\bla \ni \bla$ defined in the following way: for any box $\sq = (i, j)$ in the $k$-th partition $\lambda^k$ of $\bla$, define the vector
\[ 
    v_\sq := Y_k^j X_k^i A_k w_k \in V_{\col \sq}
\]
(We may also write $v_\sq^\bla$ when we want to stress that $\sq$ is a box in $\bla$.)
Let $U_\bla \subset \CM_\bd$ be the open subset where $\{v_\sq\}_{\sq \in \bla}$ form a basis of $V_1, \dots, V_n$. 
Writing matrices $(X, Y, A, B)$ in the basis $\{v_\sq\}_{\sq \in \bla}$, we see that $U_\bla$ is a closed subset inside an affine space, cut out by the moment map condition $\nu = 0$. The fixed point $\bla$ is the origin of this affine space. We will refer to this as the standard coordinate on $U_\bla$.

\subsubsection{Tautological bundle}
By construction, affine Laumon spaces come with tautological bundles $\CV_i$, $i = 1, \dots, n$, defined by associating the vector spaces $V_i$ to each point in $\CM_\bd $. 

The K-theory class of the tangent bundle can be described in terms of $\CV_i$'s, see formula (3.23) of \cite{neguct2018affine}:
\begin{equation}
    \label{tangentbundle}
    [T \CM_\bd] = \sum_{i=1}^{n} \left[
        \left(1 - \frac{1}{q^2}\right)
        \left(\frac{\CV_i}{\CV_{i-1}} - \frac{\CV_i}{\CV_i}\right)
        + \frac{\CV_i}{u_i^2}
        + \frac{u_{i+1}^2}{\CV_i q^2}
    \right]
\end{equation}
where $\CV_{n+i} = \CV_i p^{-2}$, $u_{n+i} = u_i p^{-1}$ for $i \in \BZ$. 

Given a partition $\bla = (\lambda^1, \dots, \lambda^n)$, for a box $\sq$ in the $(x, y)$ position of $\lambda^k$, 
let $i$ be the color of the box (i.e. $y+k \equiv i$ mod $n$),
define the weight of the box to be 
\[ 
    \chi_\sq := u_k^2 q^{2x} p^{2  \frac{y+k-i}{n} }
\]
\hide{Note that we have defined color of the box before. A box get $p^2$ in weight when going from color $n$ to color 1.}
With this notation, we have 
\begin{equation}\label{eqn:V_sum_box}
    \CV_i|_\bla = \sum_{\sq \in \bla}^{\col \sq = i} \chi_\sq
\end{equation}

As an example, the partition in Figure \ref{fig:multipartition_eg} corresponds to $n = 3$, $\bd = (4, 5, 5)$. Boxes of colors 1, 2, and 3 are colored as pink, green, and yellow, respectively.
The weight $\chi_\sq$ of each box is labeled in the box.

\ytableausetup{centertableaux, boxsize=1.8em}
\begin{figure}[htbp]
    \centering
    {\Large
        \begin{ytableau}
            *(pink) u_1^2 p^2 \\
            *(yellow) u_1^2 \\
            *(green) u_1^2 & *(green) u_1^2 q^2  \\
            *(pink) u_1^2 & *(pink) u_1^2 q^2 
        \end{ytableau}
        \: \:\:\: \: \:\:\:
        \begin{ytableau}
            \none \\
            \none \\
            *(yellow) u_2^2 \\
            *(green) u_2^2 & *(green) u_2^2 q^2
        \end{ytableau}
        \: \:\:\: \: \:\:\:
        \begin{ytableau}
            *(yellow) u_3^2 p^2 \\
            *(green) u_3^2 p^2 \\
            *(pink) u_3^2 p^2 \\
            *(yellow) u_3^2 & *(yellow) u_3^2 q^2
        \end{ytableau}   
    }
\caption{}
\label{fig:multipartition_eg}
\end{figure}

\subsection{Quantum affine algebra \texorpdfstring{$\glnhat$}{glnhat}} 

\subsubsection{Generators and relations of \texorpdfstring{$\glnhat$}{glnhat}}
\label{sec:gen_rel_uq}
We recall some facts about the quantum affine algebra $\glnhat$.
Our exposition follows Section 2.11 of \cite{neguct2013quantum}, 

One way to define $\glnhat$ is to use root generators $e_{[i;j)}$ and $f_{[i;j)}$.
For $1 \leq i \leq n$ and $j > i$, let $[i;j)$ denote the length-$n$ vector whose $k$-th entry
is equal to the number of integers congruent to $k \mod n$ in the interval $[i, j)$. 
(For example, when $n=3$, $[1;5) = (2,1,1)$). Define the positive half of $\glnhat$ by
\begin{equation}
    \uup := \langle e_{[i;j)} \rangle_{i < j}^{1 \leq i \leq n} / \text{relations}
\end{equation}
where the relations are 
\begin{multline} \label{eij_relation}
    q^{c_1(i,j,i',j')} e_{[i;j)} e_{[i';j')} - q^{c_2(i,j,i',j')} e_{[i';j')} e_{[i;j)} \\= (q - q^{-1})\left(\sum_{i' < a \leq j'}^{a \equiv i}
    e_{[a;j')} e_{[i+i'-a;j)}
    - \sum_{i' \leq a < j'}^{a \equiv j} e_{[i;j+j'-a)} e_{[i';a)}
    \right)
\end{multline}
where $c_1(i,j,i',j')$ and $c_2(i,j,i',j')$ are integers determined by $i,j,i',j'$ mod $n$.
(See formula (2.35), (2.36) of \loccit for explicit expressions of $c_1$ and $c_2$.)
Similarly, $\uum$ is defined with generators $\{f_{[i;j)}\}_{i < j}^{1 \leq i \leq n}$ with analogous relations.
These generators are natural in the geometric setting -- we will see below that 
they can be realized using correspondences on affine Laumon spaces. 

The quantum group $\glnhat$ admits a triangular decomposition
\begin{equation}
    \glnhat = \uup \otimes \BQ(q)[\psi_1^\pm, \dots, \psi_n^\pm, c^\pm] \otimes \uum
\end{equation}
See formulas (2.33), (2.34), and (2.46) in \cite{neguct2013quantum} for the relations among the tensor factors.

We will also need another set of generators of $\uup$ and $\uum$.
It was proved in \cite{neguct2013quantum} Proposition 2.21 that  
\[ 
    \glnhat \simeq \slnhat \otimes U_q(\widehat{\fgl}_1).
\]
$e_{[i;i+1)}, i=1,2,\dots,n$ generate $U_q^{+}(\widehat{\fsl}_n)$. 
Let $e_i := e_{[i;i+1)}$ and let $P_k, k=1,2,\dots$ be the generators of $U_q^{+}(\widehat{\fgl}_1)$. Then $e_i, i=1,2,...,n$ and $P_k, k=1,2,\dots$ generate $\uup$. 
Similarly, for the negative half, $f_{[i;i+1)}, i=1,2,\dots,n$ generate $U_q^{-}(\widehat{\fsl}_n)$. Let $f_i := f_{[i;i+1)}$ and let $P_{-k}, k=1,2,\dots$ be the generators of $U_q^{-}(\widehat{\fgl}_1)$.
Then $f_i, i=1,2,...,n$ and $P_{-k}, k=1,2,\dots$ generate $\uum$.

\subsubsection{Anti-involution induced from the shuffle algebra}
We will mainly consider the action of $\glnhat$ in this paper, but the shuffle algebra
is needed to define the anti-involution of $\glnhat$ and to examine the denominators
that appear in the $\glnhat$ action. So we recall some basics facts about it.
See \cite{neguct2013quantum} for details.

For a given integer $n$, the shuffle algebra $\CA$ (depending on $n$) is a Hopf algebra over $\BQ(q, p)$. It admits a triangular decomposition
\[ 
    \CA = \CA^+ \otimes \CA^0 \otimes \CA^-.
\]
The positive half $\CA^+$ has a basis consisting of color symmetric rational functions satisfying certain conditions. And the negative half $\CA^-$ has basis labeled by the same rational functions.
$\CA$ has an anti-involution $\tau$
which sends a function in the positive/negative half to the same polynomial in 
the other half.

Consider the slope 0 subalgebra
\[ 
    \glnhat \hookrightarrow \CA,
\]
as defined in \cite{neguct2013quantum}
(also called horizontal subalgebra when $\CA$ is identified with the quantum toroidal algebra).
The anti-involution $\tau $ preserves 
$\glnhat \otimes_{\BQ(q)} \BQ(q, p)$ and thus induces an anti-involution on it (still denoted $\tau$).

Note that although $\glnhat$ has nothing to do with $p$, the anti-involution 
$\tau$ does. This is because
\[ 
    \glnhat = \slnhat \otimes U_q(\widehat{\fgl}_1),
\]
and the anti-involution on $\slnhat$ and on $U_q(\widehat{\fgl}_1)$ can be scaled independently, giving
a family of different anti-involutions of $\glnhat$.
Setting $p$ to a number corresponds to choosing one anti-involution in this family.
The $\tau$ defined here is natural in this setting because the associated Shapovalov form 
can be identified with the Euler character pairing of K-theory twisted by a line bundle, see Section \ref{action_on_k}.

\subsubsection{Universal Verma module and Shapovalov form}
Recall the definition of $\BZ_\bu, \BF_\bu$ in equations (\ref{Zu}) and (\ref{Fu}).
Let
\begin{equation}
    \label{univ_verma}
    V^\text{loc} := \glnhat \otimes_{U_q^{\leq 0}(\widehat{\fgl}_n)} \BF_\bu
\end{equation}
be the universal Verma module of $\glnhat$, where the negative half $\uum$ acts on $\BF_\bu$ 
by 0, $\psi_i$ acts by scalars $u_i / q^i$, and $c$ acts by the scalar $q^n p$. The superscript ``loc'' stresses
that it is defined over the fraction field $\BF_\bu$.

The Shapovalov form on $V^\text{loc}$ is defined by 
\begin{align} \label{shapo_on_V}
    (|\emptyset \rangle, |\emptyset \rangle) &= 1, \\
    (g |u \rangle, |v \rangle) &= (|u \rangle, \tau(g) |v \rangle), \nonumber
\end{align}
where $\emptyset$ stands for the lowest weight vector in $V^\text{loc}$.

\subsubsection{PBW basis} \label{pbw_int_form}
We define a PBW basis for $\uup$.
To each multi-partition $\bla$, we define an element $e_\bla \in \uup$ as follows:
split $\bla$ into vertical strips, then multiply the root generators corresponding to 
each strip from left to right. 

\begin{example}
Let $\bla$ be the partition given in Figure \ref{fig:multipartition_eg}.
The first vertical strip has 4 boxes, with colors from 1 to 4, so it corresponds to the root generator $e_{[1;5)}$. The second vertical strip has 2 boxes, with colors 1 and 2, so it corresponds to the root generator $e_{[1;3)}$. Going through also columns, we see that  
\[ 
    e_\bla = e_{[1;5)} e_{[1;3)} e_{[2;4)} e_{[2;3)} e_{[3;7)} e_{[3;4)}.
\]
\end{example}

By Section 2.16 of \cite{neguct2013quantum},
for each $\bd \in \BZ_{\geq 0}^n$, the set $\{e_\bla\}_{|\bla| = \bd}$ forms a basis of the 
degree $\bd$ piece of $\uup$, and any product of $e_{[i;j)}$'s can be written as a 
$\BZ[q, q^{-1}]$ combination of $\{e_\bla\}$.

Let $|e_{\bla}\rangle := e_\bla |\emptyset \rangle \in V^\text{loc}$. We will call these the PBW basis of the Verma module. Let $|e_{\bla}\rangle^* \in V^\text{loc}$ denote their dual under the Shapovalov form, i.e.,
\[ 
    (|e_{\bla}\rangle, |e_{\bmu}\rangle^*) = \delta_{\bla \bmu}.
\]

\subsection{Action on \texorpdfstring{$K_\sT(\CM_\bd )$}{ktmd}} \label{action_on_k}
As shown in \cite{braverman2005finite,neguct2018affine,tsymbaliuk2010quantum}, 
there is a natural action of $\glnhat$ on $K^\text{loc}$ defined by correspondences. We briefly recall the construction and refer to \cite{neguct2018affine}, sections 3 and 4, for more details.

\begin{remark}\label{larger_alg}
In fact, it is shown in \loccit that there is an action of quantum toroidal algebra $U_q(\widehat{\widehat{\fgl}}_n)$ on $K ^\text{loc}$ by tensoring with tautological line bundles. In the present paper, we will only use the subalgebra $\glnhat$. Similarly, for finite Laumon space $\CM_\bd ^{\text{fin}}$, there is an action of  $U_q(\widehat{\fsl}_n)$, but we will only use the subalgebra $U_q(\fsl_n)$.
\end{remark}

The geometric realization relies on certain varieties 
\footnote{To be precise, the correspondences $\fZ_{[i;j)}$, $\ofZ_{[i;j)}$ both have 
connected components labeled by $\bd \in \BZ _{\geq 0}^{n-1}$ so that $\fZ_{[i;j)}^{\bd }$ and $\ofZ_{[i;j)}^{\bd }$
maps to $\CM_{\bd + [i;j)}$ and $\CM_{\bd}$. We omit the superscript $\bd $ to declutter the notation.}
$\fZ_{[i;j)}$, $\ofZ_{[i;j)}$ together with 
K-theory classes $[\fZ_{[i;j)}^+] \in K_\sT(\fZ_{[i;j)})$, $[\ofZ_{[i;j)}^-] \in K_\sT(\ofZ_{[i;j)})$ and maps
\begin{equation*}
    \xymatrix{ & \fZ_{[i;j)} \text{ or } \ofZ_{[i;j)} \ar[ld]_{p^+ \text{ or } \op^+} \ar[rd]^{p^- \text{ or } \op^-}  & \\
    \CM_{\bd^+} & & \CM_{\bd^-}}
\end{equation*}
where $\bd ^+ - \bd ^- = [i;j)$. See \cite{neguct2018affine}, formula (4.13),(4.16)
for the exact definitions of the classes $[\fZ_{[i;j)}^+]$, $[\ofZ_{[i;j)}^-]$.

The action of the root generators are given by 
\begin{align*}
    e_{[i;j)}(\alpha) &= p^+_*([\fZ_{[i;j)}^+] \cdot p^{-*}(\alpha)), \\
    f_{[i;j)}(\alpha) &= \op^-_*([\ofZ_{[i;j)}^-] \cdot p^{+*}(\alpha)),
\end{align*}
for any class $\alpha \in K^\text{loc}$. 
Also define the Cartan part action by  
\begin{equation} \label{eqn:psi_action}
    \psi _i = \text{ multiplication by }\frac{q ^{i + d _{i-1}-d_i}}{u_i} \text{ on } K_\bd ^\text{loc},
\end{equation} 
\begin{equation}
    c = \text{ multiplication by } q^n p.
\end{equation}
These operators generate the action of $\glnhat$ on $K^\text{loc}$.

The action can also be written explicitly in the normalized fixed-point basis.
The action of $e_{[i;j)}$ is given by
\begin{equation}
    \label{actioncoeff}
    \langle \bla | e_{[i;j)} | \bmu \rangle = R^+_{ij}(\bla \backslash \bmu) \prod_{\bsq \in \bla \backslash \bmu} 
    \left[(q^{-1}-q) \zeta\left(\frac{\chi_\bsq}{\chi_\bmu}\right) \tau_+(\chi_\bsq)\right],
\end{equation}
where   
\begin{equation} \label{rij}
    R^+_{ij} = \text{Sym} \left[\frac{1}{(1 - \frac{z_{i+1}}{z_i q^2}) \cdots (1 - \frac{z_{j-1}}{z_{j-2} q^2})}
    \prod_{i \leq a < b < j} \zeta\left(\frac{z_b}{z_a}\right)\right].
\end{equation}
See \cite{neguct2018affine} for the exact definitions of $\zeta$ and $\tau_+$ and how to evaluate $R_{ij}$ on $\bla \backslash \bmu$. The formula above for $e_{[i;j)}$ is formula (3.59) in \cite{neguct2018affine} and the formula for $R^+_{ij}$ can be deduced from formula (2.12)-(2.16) there.
In particular, the product of $\zeta$ and $\tau_+$ above can be simplified (\cite{neguct2018affine}, formula (3.56)): 
\begin{equation}
    \zeta\left(\frac{z}{\chi_\bla}\right) \tau_+(z) = 
    \frac{\prod^{\col \sq = \col z + 1}_{\text{addable } \sq \text{ of } \bla} 
    \left(\frac{\sqrt{\chi_\sq}}{q} - \frac{zq}{\sqrt{\chi_\sq}}\right)}{\prod^{\col \sq = \col z}_{\text{removable } \sq \text{ of } \bla} 
    \left(\sqrt{\chi_\sq} - \frac{z}{\sqrt{\chi_\sq}}\right)}. \label{simple_zeta_tau}
\end{equation}
A similar formula holds for the action of $f_{[i;j)}$.

\begin{remark}\label{rmk:box_color}
Note that $\langle \bla | e_{[i;j)} | \bmu \rangle \neq 0$ implies that $\bmu$ is 
contained in $\bla$, and the difference $\bla \backslash \bmu$ has boxes with colors given by $[i;j)$. This property will be used frequently in the rest of the paper.
\end{remark}


\begin{lemma}
    \label{corr_proper}
    The maps $p^+, p^-, \op^+$ are proper. 
\end{lemma}
\begin{proof}
For $p^+$ and $p^-$, \cite{neguct2018affine} showed that, in the moduli-of-sheaves description of affine Laumon spaces $\CM_\bd $, 
$\fZ _{[i;j)}$ can be realized as a
chain of modifications on the flag of sheaves. The usual argument for Hecke correspondences then shows that $p^+$ and $p^-$ are proper.

To show the properness of $\op^+$, we first recall the definition of $\ofZ_{[i;j)}$ in Section 4.6 of \cite{neguct2018affine}. 
Fix $\bd^+, \bd^-$ such that $\bd^+ - \bd^- = [i;j)$. Fix two sets of vector spaces $\bV^+, \bV^-$ of dimensions $\bd^+$, $\bd^-$ and a full
flag of subspaces
\begin{equation}
\label{eqn:zeit}
0 = U^{[i;i)} \subset U^{[i;i+1)} \subset U^{[i;i+2)} \subset ... \subset U^{[i;j)} = \text{Ker}(\bV^+ \twoheadrightarrow \bV^-)
\end{equation}
where each $U ^{[i;a)}$ is an $n$-tuple of vector spaces of dimensions $[i;a)$.
View $M_{\bd^+}$ (defined in (\ref{prequotient})) as linear maps $(\bX, \bY, \bA, \bB)$ acting on $\bV^+$, and consider the subspace $Z_{[i;j)} \subset M_{\bd^+}$ where $\bX, \bY$ satisfy the condition that 
\begin{equation}\label{eqn:cl_cond}
    \bX(U^{[i;a)}) \subset U^{[i;a-1)}, \quad \bY(U^{[i;a)}) \subset U^{[i;a+1)} \text{ for } a = i+1, \dots, j-1.
\end{equation}
Let $\eta$ be the moment map on $Z_{[i;j)}$ induced from $\nu$ on $M_{\bd^+}$. Let $P_{[i;j)}$ be the parabolic subgroup that preserves the flag of subspaces (\ref{eqn:zeit}). Define
\[ 
    \ofZ_{[i;j)} = \eta^{-1}(0)^s / P_{[i;j)}.
\] 
Since the conditions in \eqref{eqn:cl_cond} are closed conditions, $\eta^{-1}(0)^s$ is a closed subvariety of $\nu^{-1}(0)^s$. Thus,
$\ofZ_{[i;j)}$ is a closed subvariety of $\nu^{-1}(0)^s / P_{[i;j)}$, which is a projective bundle over $\CM_{\bd^+} = \nu^{-1}(0)^s / G_{\bd^+}$. Therefore, the projection map $\op^+$ is proper.

\end{proof}

\begin{remark} \label{rmk:tau_proper}
Note that (see \cite{neguct2018affine} formula (4.13)-(4.14)) both $e _{[i;j)}$
and its anti-involution $\tau(e_{[i;j)})$ are defined using K-theory classes on $\fZ_{[i;j)}$, so both map $K ^{\text{int}}$ to $K ^{\text{int}}$.
\end{remark}

\begin{remark}
    The map $\op^-$ is not proper in general. This will be further discussed in section \ref{action_almost_preserve}.
\end{remark}

There is an isomorphism (\cite{neguct2018affine}, see also \cite{braverman2005finite}, \cite{tsymbaliuk2010quantum}) of $\glnhat$ modules
\begin{equation}
    \label{verma_iso}
    \Phi : V^\text{loc} \xrightarrow{\simeq} K^\text{loc}
\end{equation}
sending the lowest weight vector in $V^\text{loc}$ to the structure sheaf in degree 0.
In the rest of this paper, we will use $|e_{\bla}\rangle, |e_{\bla}\rangle^*$ to denote either 
an element in $V^\text{loc}$ or its image in $K^\text{loc}$ under the isomorphism $\Phi$ (and it should be clear which one it is).

The Shapovalov form on $V^\text{loc}$ induces a bilinear form on $K^\text{loc}$ given by (cf. \cite{braverman2005finite}, Proposition 2.29)
\footnote{
The reason that we use $\mathcal{D}_\bd^{-1}$ here as opposed to $\mathcal{D}_\bd$ used in \cite{braverman2005finite} is that the correspondences we use (see formula (4.13)-(4.16) of \cite{neguct2018affine}, especially the $j-i=1$ case) are slightly different from the ones in \cite{braverman2005finite} Section 2.8 in the twists by the tautological line bundles $\CL_i$. 
}:
\begin{equation}\label{shapo_k}
    (\CF, \CG) = \const_d \cdot \chi(\CF \otimes \CG \otimes \mathcal{D}_\bd^{-1}),
\end{equation}
where $\chi $ is the (equivariant) Euler characteristic and
\[ 
    \mathcal{D}_\bd = \prod_{i=1}^n \det \CV_i
\]
is the determinant line bundle of tautological bundles, and $\const_d$ is an invertible element in $\BZ_\bu$ (i.e., a monomial with coefficient 1) depending on $d$.

\subsection{Specializing to finite Laumon spaces and \texorpdfstring{$U_q(\fsl_n)$}{uqsln} }\label{fin_laumon}
When $d_n = 0$, the affine Laumon space becomes the finite Laumon space 
initially defined in \cite{laumon1988analogue}, 
which we denote by $\CM^\text{fin}_\bd$, where $\bd = (d_1, \dots, d_{n-1})$, see \cite{neguct2018affine}, Remark 3.3. 
In this section, we spell out what the above constructions mean in this special case.

The quiver becomes (cf. \cite{nakajima2011handsaw}, Section 2(i)):

\begin{picture}(200,110)(40,-30)\label{handsaw}

    \put(115,31){\vector(1,0){50}}
    \put(135,34){$Y_{2}$}
    
    \put(175,31){\vector(1,0){40}} 
    \put(195,34){$Y_{3}$}
    
    \put(247,31){\vector(1,0){38}}
    \put(255,34){$Y_{n-1}$}
    
    \put(110,31){\circle*{10}}
    \put(110,50){\circle{30}}
    \put(95,47){$X_{1}$}
    \put(117,22){$\color{red}{V_{1}}$}
    \put(102,36){\vector(4,-1){5}}
    
    \put(170,31){\circle*{10}}
    \put(170,50){\circle{30}}
    \put(155,47){$X_{2}$}
    \put(177,22){$\color{red}{V_{2}}$}
    \put(162,36){\vector(4,-1){5}}
    
    \put(225,23){${\cdots}$}

    \put(290,31){\circle*{10}}
    \put(290,50){\circle{30}}
    \put(275,47){$X_{n-1}$}
    \put(297,22){$\color{red}{V_{n-1}}$}
    \put(282,36){\vector(4,-1){5}}
    
    \put(75,-20){\line(1,0){10}}
    \put(75,-10){\line(1,0){10}}
    \put(75,-20){\line(0,1){10}}
    \put(85,-20){\line(0,1){10}}
    \put(88,-20){$\color{blue}{\BC w_{1}}$}
    
    \put(85,-10){\vector(2,3){25}}
    \put(98,4){$A_{1}$}
    
    \put(135,-20){\line(1,0){10}}
    \put(135,-10){\line(1,0){10}}
    \put(135,-20){\line(0,1){10}}
    \put(145,-20){\line(0,1){10}}
    \put(148,-20){$\color{blue}{\BC w_{2}}$}
    
    \put(111,27){\vector(2,-3){25}}
    \put(145,-10){\vector(2,3){25}}
    \put(128,4){$B_{2}$}
    \put(158,4){$A_{2}$}
    
    \put(195,-20){\line(1,0){10}}
    \put(195,-10){\line(1,0){10}}
    \put(195,-20){\line(0,1){10}}
    \put(205,-20){\line(0,1){10}}
    \put(208,-20){$\color{blue}{\BC w_{3}}$}
    
    \put(171,27){\vector(2,-3){25}}
    \put(205,-10){\vector(2,3){15}} 
    \put(188,4){$B_{3}$}
    \put(216,2){$A_{3}$}
    
    \put(255,-20){\line(1,0){10}}
    \put(255,-10){\line(1,0){10}}
    \put(255,-20){\line(0,1){10}}
    \put(265,-20){\line(0,1){10}}
    \put(268,-20){$\color{blue}{\BC w_{n-1}}$}
    
    \put(239,12){\vector(2,-3){15}}
    \put(265,-10){\vector(2,3){25}}
    \put(248,4){$B_{n-1}$} 
    \put(278,4){$A_{n-1}$}
    
    \put(315,-20){\line(1,0){10}}
    \put(315,-10){\line(1,0){10}}
    \put(315,-20){\line(0,1){10}}
    \put(325,-20){\line(0,1){10}}
    \put(328,-20){$\color{blue}{\BC w_{n}}$}
    
    \put(291,27){\vector(2,-3){25}}
    \put(308,4){$B_{n}$}
    
    
\end{picture}

\noindent Recall the description of fixed points of affine Laumon spaces in section \ref{fixedpt}.
Specializing to the finite case is to require that there is no box of color $n$ in the partition.
In other words, the fixed points in the finite case are parametrized by $n-1$ tuples of 2d partitions $\bla = (\lambda^1, \dots, \lambda^{n-1})$, where each $\lambda^i$ is a 2d partition such that $\lambda^i$ has at most $n-i$ rows. We will call such $\bla$ a Kostant partition, following the nomenclature in \cite{finkelberg2000parabolic}. 
Figure \ref{fig:partition_finite} shows an example of $n=4$ with degree $\bd = (4, 4, 3)$. Boxes of color 1, 2, and 3 are colored as pink, green, and yellow, respectively.
\ytableausetup{centertableaux, boxsize=2em}
\begin{figure}[htbp]
    \centering
    {\Large
        \begin{ytableau} 
            *(yellow) u_1^2 \\
            *(green) u_1^2 & *(green) u_1^2 q^2  \\
            *(pink) u_1^2 & *(pink) u_1^2 q^2 & *(pink) u_1^2 q^4 & *(pink) u_1^2 q^6 
        \end{ytableau}
        \: \:\:\:
        \begin{ytableau}
            \none    \\
            *(yellow) u_2^2 \\
            *(green) u_2^2 &  *(green) u_2^2 q^2
        \end{ytableau}
        \: \:\:\:
        \begin{ytableau}
            \none   \\
            \none   \\
            *(yellow)  u_3^2
        \end{ytableau}
    }   
\caption{}
\label{fig:partition_finite}
\end{figure}

For the quantum group action, the root generators $e_{[i;j)}, f_{[i;j)}$ with $1 \leq i < j \leq n$ preserve $K_\sT(\CM_\bd^\text{fin})$ since they will not add or remove color $n$ boxes (c.f. Remark \ref{rmk:box_color}). Together with $\psi _{i+1}/ \psi_i, i=1,2,...,n-1$, they generate a subalgebra
\[ 
    U_q(\fsl_n) \hookrightarrow \glnhat.
\]
The generators $e_{[i;i+1)}$ correspond to the usual generators of $U_q(\fsl_n)$ multiplied by $q - q^{-1}$. From relation (\ref{eij_relation}), we have
\[ 
    q e_{[i+1;j)} e_{[i;i+1)} - e_{[i;i+1)} e_{[i+1;j)} = (q - q^{-1}) e_{[i;j)}
\]
for $1 \leq i < j \leq n$.

In what follows, we will also use that for finite Laumon spaces, the relative normal bundles can be expressed as 
\footnote{Note that in the finite case we have $$\left(1-\frac{1}{q^2}\right)\sum _{i \leq a < b <j}^{a \equiv b-1} \frac{\CL_a}{\CL_b} + \sum_{a=i+1}^{j-1} \frac{\CL_a}{\CL_{a-1}q^2}=\sum_{a=i+1}^{j-1} \frac{\CL_a}{\CL_{a-1}}.$$ }
(c.f. \cite{neguct2018affine}, Proposition 4.14, 4.15)
\begin{multline} \label{relweights}
[T\fZ_{[i;j)}] - p ^{+*} [T\CM_{\bd^+} ] = i-j + \\
\left(1 - \frac{1}{q^2}\right) \left(\sum_{i \leq a < j} \frac{\CV_a^+}{\CL_a}
- \sum_{i \leq a < j}  \frac{\CV_{a+1}^+}{\CL_a}\right)
+ \sum_{a=i+1}^{j-1} \frac{\CL_a}{\CL_{a-1}}
- \sum_{a=i}^{j-1} \frac{u_{a+1}^2}{\CL_a q^2}.
\end{multline}
\begin{multline} \label{relweights2}
[T\fZ_{[i;j)}] - p ^{-*}[T\CM_{\bd^-} ] = i-j + \\
\left(1 - \frac{1}{q^2}\right) \left(-\sum_{i \leq a < j} \frac{\CL_{a}}{\CV_a^-} 
+ \sum_{i \leq a < j} \frac{\CL_a}{\CV_{a-1}^-}\right)
+ \sum_{a=i+1}^{j-1} \frac{\CL_a}{\CL_{a-1}}
+ \sum_{a=i}^{j-1} \frac{\CL_a}{u_a^2}.
\end{multline}
\hide{$\fZ$ and $\ofZ$ are the same in finite case. The reason that $f_{[i;j)} \neq e_{[i;j)}$ when $j-i>1$
is that the twisting by line bundles in the fundamental class.}
\hide{Also, we have simplified $(1-q^{-2})\sum {L_a/L_b} + \sum {L_a/L_{a-1}q^2}$ to $\sum {L_a/L_{a-1}}$}
Here, $\CV_a^+$ (resp. $\CV_a^-$) is the tautological bundle for $\CM_{\bd^+}$ (resp. $\CM_{\bd^-}$), and $\CL_a$ is their difference on $\fZ_{[i;j)}$. (There are also analogous formulas for $\ofZ_{[i;j)}$, but we won't use them in this paper.)

\section{\texorpdfstring{$K_\sT(\CM_\bd^\text{fin} )$}{ktmdfin} and dual Verma module of \texorpdfstring{$U_q(\fsl_n)$}{uqsln}} \label{proof_finite_case}

In this section, we will focus on the equivariant K-theory of finite Laumon space $K_\sT(\CM_\bd^\text{fin} )$, and all notations will be adapted to the finite case: $\bd$ will denote a vector of length $n-1$; boldface letters
(e.g., $\bla$) will denote Kostant partitions as explained in section \ref{fin_laumon};
$K ^\text{loc}, K ^\text{int}, V ^\text{loc}, \Phi$ are understood as defined for 
$\CM_\bd ^{\text{fin}}$ and $U_q(\fsl_n)$. (Note that $\Phi $ is still an isomorphism, see Section 2.26 of \cite{braverman2005finite}.)


\begin{theorem}
    \label{maintheorem}
    (1) For any partition $\bla$, $\Phi$ maps the dual PBW basis $|e_{\bla}\rangle^*$ into 
    $K_\sT(\CM_\bd^{\text{fin}} )$ where $\bd =|\bla|$. 

    (2) For any $\bd \in \BZ_{\geq 0}^{n-1}$, $\{\Phi \left(|e_{\bla}\rangle^*\right)\}_{|\bla|=\bd }$ form a $\BZ_\bu$-basis
    of $K_\sT(\CM_\bd ^{\text{fin}})$.
\end{theorem}

\begin{remark}
A similar result was obtained in \cite{feigin2011gelfand}, Theorem 3.5 when studying the equivariant cohomology of the Laumon space as a (non-quantum) $U(\mathfrak{gl}_n)$-module. The approach we use here is different from theirs.
\end{remark}

Hence, specializing equivariant parameters, we can get dual Verma module of a given lowest weight. 
The lowest weight is determined by the action of $\psi _{i+1}/\psi_i, i=1,...,n-1$, see \eqref{eqn:psi_action}. So we have 
\begin{corollary}
    \label{dual_verma_corr}
    Fix $a_1,...,a_n \in \BC$ and let $q_0 \in \BC$ be generic. 
    Let $\BC_a$ denote the evaluation module of $\BZ_\bu$, i.e., the target of the map
    \[ 
        \BZ_\bu=\BZ[u_1^\pm,...,u_n^\pm, q^\pm] \xrightarrow{u_i = q ^{a_i},q=q_0} \BC
    \]
    Then $K ^\text{int}\otimes _{\BZ_\bu} \BC_a$ is isomorphic to the dual Verma module of $U_{q_0}(\fsl_n)$ with lowest weight $(b_1,...,b_{n-1})$ where 
    $b_i = a_i - a _{i+1}+1$.
\end{corollary}

\subsubsection*{Proof of Theorem \ref{maintheorem}}
    
First, observe that part (1) of the theorem implies part (2).
Indeed, given any $\alpha \in K_\sT(\CM_\bd^\text{fin} )$ and any PBW basis $|e_{\bla}\rangle$, using the pairing (\ref{shapo_k}), we have:
\begin{equation*}
        (\alpha , \Phi (|e_{\bla}\rangle)) = (\tau (e_\bla)\alpha, \Phi (|\emptyset \rangle ) ).
\end{equation*}
Note that $\tau (e_\bla)\alpha$ is an integral K-theory class, see Remark \ref{rmk:tau_proper}. Thus, the pairing lives in $\BZ_\bu$. By definition, the dual PBW basis satisfies
\[ 
    (\Phi (|e_{\bla}\rangle^*), \Phi (|e_{\bmu}\rangle)) = \delta _\bla^\bmu.
\]
In other words, for any $\alpha \in K_\sT(\CM_\bd ^{\text{fin}})$,
\[ 
    \alpha = \sum_{|\bla|=\bd }(\alpha , \Phi (|e_{\bla}\rangle))\cdot \Phi (|e_{\bla}\rangle^*).
\]
So $\alpha$ is in the $\BZ_\bu$-linear span of $\{\Phi (|e_{\bla}\rangle^*)\}_{|\bla|=\bd }$. 

Thus, we only need to prove that $\Phi (|e_{\bla}\rangle^*) \in K_\sT(\CM_\bd^\text{fin} )$. 
For any partition $\bla$, let $h(\bla)$ denote the sum of the heights of all boxes in $\bla$, where the height of a box is defined to be its row index (starting from 0). Choose a total ordering of the fixed points such that if $h(\bmu) < h(\bla)$ then $\bmu < \bla$. (If $h(\bmu)=h(\bla)$ then they can be in either order.)
In this order, the flat partition (i.e., each partition has only one row) is the smallest. 

\begin{lemma}
    \label{pbw_vanish}
    Assume that $|e_{\bla}\rangle|_\bmu \neq 0$ (i.e. the restriction of $\Phi (|e_{\bla}\rangle)$ to the fixed point $\bmu$ is non-zero). Then $h(\bmu) \leq h(\bla)$, and $h(\bmu) = h(\bla)$ only when $\bmu = \bla$.
\end{lemma}
\begin{proof}
    Recall that formula (\ref{actioncoeff}) implies that $ \langle \bla|e_{[i;j)}|\bmu \rangle $ is only non-zero if $\bla$ can be obtained from $\bmu$ by adding boxes of color $i, i+1,...,j-1$. 
    
    We prove the lemma by induction. Assume the lemma is proven for all partitions of degree $\bd '$ with $|\bd '| < |\bd |$ (here $|\bd |=d_1 + ... + d _{n-1}$), 
    we show that it holds for partitions of degree $\bd$.

    For a partition $\bla$, assume that the left most vertical strip in $\bla$ contains boxes of color from $i$ to $j-1$. 
    (This implies that the first non-empty partition in $\bla$ is $\bla^i$.)
    Let $\bla'$ be $\bla$ with the leftmost vertical strip $[i;j)$ removed. In terms of PBW basis, $|e_{\bla}\rangle = e _{[i;j)}|e_{\bla'}\rangle$. 
    Then
    \[ 
        h(\bla) = h(\bla') + 0 + 1 + ... + (j-i-1).
    \]
    By the induction hypothesis, if $|e_{\bla'}\rangle|_{\bmu'} \neq 0$ for some $\bmu'$, 
    then $h(\bmu') \leq h(\bla')$ and $h(\bmu') = h(\bla')$ only when $\bmu' = \bla'$. 
    We examine what happens when we apply $e_{[i;j)}$ to a fixed point class $|\bmu' \rangle $.
    If $e_{[i;j)} |\bmu' \rangle|_\bnu \neq 0 $, then $\bnu$ is obtained from $\bmu'$ by adding boxes of color $i, i+1,...,j-1$. The maximal height that $\bnu$ could have is when adding these boxes all in the $i$-th partition, which implies 
    \[ 
        h(\bnu) = h(\bmu') + 0 + 1 + ... + (j-i-1).
    \]
    So $h(\bnu) \leq h(\bla)$. The equality holds only when $\bmu' = \bla'$ and all boxes are added to the $i$-th partition, i.e., $\bnu = \bla$.
    
\end{proof}

Dual PBW basis is uniquely determined by the property that  
\[ 
    (|e_{\bla}\rangle^*, |e_{\bmu}\rangle) = \delta _{\bla}^{\bmu}
\]
For each partition $\bla$ with $|\bla|=\bd$, we will construct elements $g _{\bla} \in K_\bd ^\text{int}$ 
satisfying the following properties:
\begin{align}
    \label{vanishing}\tag{P1} &g_\bla|_\bmu \neq 0 \text{ only if } \bla \prec \bmu, \\
    \label{diag} \tag{P2} &g_\bla|_\bla \cdot |e_{\bla}\rangle|_\bla = \Lambda ^\bullet (T\CM_\bd^\text{fin} |_\bla)^\vee. 
\end{align}
Combining with Lemma \ref{pbw_vanish}, we see that for $\bmu < \bla$, $(g _{\bla}, |e_{\bmu}\rangle)$ must be 0 
since at least one of them vanishes on any fixed point. 
On the other hand, $(g _{\bla}, |e_{\bla}\rangle)$ only gets a non-zero contribution from the fixed point $\bla$, 
so (\ref{diag}) ensures that it is equal to 1 after scaling $g_\bla$ by an invertible element in $\BZ_\bu$. 
(This invertible element in $\BZ_\bu$ comes from the twisting by $\mathcal{D}_\bd^{-1}$ and the prefactor in formula (\ref{shapo_k}).)
To sum up, the inner product between a suitable $\BZ_\bu$-rescaling of $\{g_\bla\}$ and $\{|e_{\bmu}\rangle\}$ is upper-triangular with 1 on the diagonal.
\[
    \left(
    \begin{array}{ccccc}
    1                                    \\
      & 1             &   & \text{\huge *} \\
      &               & \ddots                \\
      & \text{\huge0} &   & 1            \\
      &               &   &   & 1
    \end{array}
    \right)
\]  
This will allow us to construct $|e_{\bla}\rangle^*$ as a $\BZ_\bu$-linear combination of $g_\bla$ and thus $|e_{\bla}\rangle^* \in K_\bd ^\text{int}$.

The desired class $g_\bla$ is constructed through attracting manifolds with respect to the torus action.
Recall that if a torus $\sT = (\mathbb{C}^*)^{k}$ acts on a variety $X$, given a fixed point $p \in X^\sT$ and a subtorus $\sigma : \mathbb{C}^*\to \sT$ such that $X ^{\mathbb{C}^*} = X ^{\sT}$, the attracting manifold of $p$ is defined to be
\[ 
    \text{Attr}_p := \{x \in X| \lim_{t \to 0}\sigma (t)\cdot x = p\}.
\]
In our situation, consider the torus $\sT'$ which is the torus in (\ref{torus}) without the $\mathbb{C}^*_p$ factor. ($\mathbb{C}^*_p$ acts trivially on $\CM_\bd ^{\text{fin}}$.)
Let $u_1,...,u_n, q$ be the weights of $\sT'$.
For a fixed point $\bla$ with $|\bla|=\bd $, let $\text{Attr}_\bla$ be the attracting locus of $\bla$ with respect to the subtorus 
\[ 
    \sigma_\bd  (t) = (t ^{2|\bd|},t ^{4|\bd|},...,t ^{2n|\bd|}, t ^{-1}) \subset \sT.
\]
Under this torus action, tangent weight $q^{-k}$ is attracting for $k>0$. Tangent weight $q^k u_j^2/u_i^2$ is attracting for $j>i$ and any $k$.
Let 
\[ 
    g _{\bla} := i_*[\CO _{\overline{\text{Attr}}_\bla} ] \in K_\sT(\CM_\bd^\text{fin} ).
\]
Here $\overline{\text{Attr}}_\bla$ is the closure of ${\text{Attr}}_\bla$, a closed subvariety of $\CM_\bd^\text{fin} $, and $i$ is the inclusion map. 
It remains to prove that $g_\bla$ satisfies properties (\ref{vanishing}) and (\ref{diag}).
Now the proof reduces to the following two propositions:

\begin{proposition}\label{prop:attr_vanish}
For fixed points $\bla, \bmu \in K_\sT(\CM_\bd^\text{fin} )$, if $g_\bla|_\bmu \neq 0$, then either $\bla=\bmu$ or $h(\bla) < h(\bmu)$.
\end{proposition}
\begin{proof}
Suppose on the contrary that $h(\bla) \geq  h(\bmu)$ and $\bla \neq \bmu$. 
This would imply that there exists a color $c$ and $k \leq c$ such that 
\[ 
    \sum_{i=k}^c \bla^i _{[c]} < \sum_{i=k}^c \bmu^i _{[c]},
\]
where $\bla^k _{[c]}$ is the number of color $c$ boxes in partition $\lambda ^k$. 
Indeed, if this is not the case, then for each color $c$, we have
\[ 
    \sum_{i=1}^c \bla^i _{[c]} \geq \sum_{i=1}^c \bmu^i _{[c]}.
\]
Since $\bla$ and $\bmu$ have the same number of color $c$ boxes, we can compare the height of boxes of color $c$ from right to left, and the height of the box in $\bla$ is always smaller or equal to the height of the corresponding box in $\bmu$.
This implies that $h(\bla) \leq h(\bmu)$, and the equality holds only when $\bla = \bmu$, contradicting our assumption.

Now let $\mathrm{C}$ be the set of color $c$ boxes in $\bmu^k,...,\bmu^c$. 
Consider the open subset 
\begin{equation}
    U_\mathrm{C} = \left\{\{v^\bmu_\sq\}_{\sq \in \mathrm{C}} \text{ are linearly independent}\right\}. \label{chart_uc}
\end{equation}
This is an open subset in $\CM_\bd^\text{fin} $ that contains $\bmu$. We will show that $U_C$ does not intersect with the attracting locus of $\bla$.

Consider the standard coordinates $\{v^\bla_\sq\}_{\sq \in \bla}$ on the open subset $U_\bla$ defined in Section \ref{local_coord}. Consider the vectors 
\[ 
    v _{ij}, i=1,2,...,c, j=1,2,...,\bla^i _{[c]}
\]
which form a basis of $V_c$ for points in $U_\bla$. Here, $v _{ij}$ is the vector corresponding to the $j$-th color $c$ box in partition $\lambda^i$.
For any point in $U_\bla \cap U _{\mathrm{C}}$, 
consider the vector 
\[ 
    v_0 := X_c v_{k, \bla^k_{[c]}}.
\]
By the construction of $U_C$, this vector is linearly independent of $\{v _{ij}\}_{i \geq k}$. So in the basis $\{v _{ij}\}$, $v_0$ has non-zero coefficient in some $v _{i'j'}$ for some $i' < k$. Under the action of $\sT$, this coefficient (as entry in the matrix $X_c$ in the standard basis) has weight $u_{i'}^2 / u_k^2$ multiplied by some power of $q$, which is repelling.
So this point is not in the attracting locus of $\bla$.  
Note that $U_\bla$ is 
invariant under $\sT$ so $\text{Attr}_\bla \subset U_\bla$. Thus $\text{Attr}_\bla \cap U _{\mathrm{C}} = \emptyset$.
This finishes the proof.
\end{proof}

\begin{proposition}
    $\Lambda ^\bullet  T ^\vee_\bla \mathrm{Attr}_\bla$ is equal to $|e_{\bla}\rangle|_\bla$ up to an invertible constant in $\BZ_\bu$.
\end{proposition}
Note that this implies property (\ref{diag}) because $g_\bla | _\bla = \Lambda ^\bullet  T ^\vee _{\bla, \text{repelling directions}}$.
\begin{proof}
We prove this by induction. As in the proof of Lemma \ref{pbw_vanish}, let $\bmu$ be the partition obtained from $\bla$ by removing the left most vertical strip $[i;j)$. So $|e_{\bla}\rangle = e_{[i;j)}|e_{\bmu}\rangle$. Assume that the proposition is proved for $\bmu$. Let's examine how the two expressions change for $\bla$.

Lemma \ref{pbw_vanish} implies that the coefficient $|e_{\bla}\rangle|_\bla$ comes only from $|e_{\bmu}\rangle|_\bmu$. So it changes by the exterior power of the relative tangent bundle (\ref{relweights}) up to invertible constants. So to prove the proposition, it suffices to prove that 
\begin{equation}\label{eqn:to_prove}
    \left([T\fZ_{[i;j)}] - p^{+*}[T\CM_{\bd^+} ]\right)\big|_{(\bla, \bmu)}
    = T _\bla \text{Attr}_\bla - T _\bmu \text{Attr}_\bmu.
\end{equation}
(Note that $\Lambda ^\bullet V$ and $\Lambda ^\bullet V ^\vee $ differ by an invertible constant, thus we get rid of the dual in the above formula.)

The righthand side of (\ref{eqn:to_prove}) is equal to 
the attracting weights in
\[ 
      p^{+*}[T\CM_{\bd^+}] -   p ^{-*}[T\CM_{\bd^-} ]
\]
So the proof reduces to the following 
\begin{lemma}
The weights in 
\[ 
    \left([T\fZ_{[i;j)}] - p^{+*}[T\CM_{\bd^+} ]\right)\big|_{(\bla, \bmu)}
\]
are all attracting, and the weights in
\[ 
    \left([T\fZ_{[i;j)}] - p^{-*}[T\CM_{\bd^-} ]\right)\big|_{(\bla, \bmu)}
\]
are all repelling.

\end{lemma}

\begin{proof}
To prove the first statement,
consider formula (\ref{relweights}). At the fixed point $(\bla, \bmu)$, we have, for $i \leq a < j$, 
\[ 
    \CV_a^+ |_{(\bla, \bmu)} = \sum_{k=i}^a \sum_{l=0} ^{\bla^k _{[a]} -1} u_k^2 q ^{2l}
\]
and
\[ 
    \CL_a |_{(\bla, \bmu)} = u_i^2 q^{2(\bla^i _{[a]} -1)}
\]
because all new boxes are added to the $i$-th partition. 
Recall that our choice of attracting direction is that $q^m u_k^2 / u_i^2$ is attracting for $k>i$ and any $m$, and $q^{-m}$ is attracting for $m>0$. 
Thus, most terms in 
\[ 
    \left(1 - \frac{1}{q^2}\right)\frac{\CV_a^+}{\CL_a}
\]
can be easily seen to be attracting. The only term that needs to be checked is the ones that only contains $q$, and we are left with
\[ 
    \left(1 - \frac{1}{q^2}\right) \sum_{l = -(\bla^i _{[a]} -1)}^0 q ^{2l} = 1 - q^{-2\bla^i _{[a]}}
\]
Similarly, among the terms in 
\[ 
    \left(1 - \frac{1}{q^2}\right)\frac{\CV_{a+1}^+}{\CL_a}
\]
the weights with only $q$ are
\[ 
    \left(1 - \frac{1}{q^2}\right) \sum_{l = -(\bla^i _{[a]} -1)}^{\bla^i _{[a+1]} - \bla^i _{[a]}} q ^{2l} = q ^{2\bla^i _{[a+1]}-2\bla^i _{[a]}} - q^{-2\bla^i _{[a]}} 
\]
Since $\bla^i _{[a+1]} \leq \bla^i _{[a]}$, both of the above terms are attracting. (There could be constant weight 1 appearing in this expression. But they will eventually cancel out because all fixed points are isolated.)
The terms in 
\[ 
    \frac{\CL_a}{\CL_{a-1}}
\]
are easily seen to be attracting. So we proved the first statement.

The proof for the second statement is similar. The only difference is that the term 
\[ 
    \frac{\CL_a}{\CL_{a-1}}
\]
will contribute attracting weights
\[ 
    q ^{2\bla^i _{[a]}-2\bla^i _{[a-1]}}
\]
But these weights will cancel with the second term from
\[ 
    \left(1 - \frac{1}{q^2}\right)\frac{\CL_a}{\CV_{a-1}^-} = q ^{2\bla^i _{[a]}-2} - q ^{2\bla^i _{[a]} - 2\bla^i _{[a-1]}}
\]
Other terms in \eqref{relweights2} are all repelling.
\end{proof}

\end{proof}

\section{$K_\sT(\CM_\bd)$ and dual Verma module of $\glnhat$} \label{proof_aff_case}

We will prove the analog of Theorem \ref{maintheorem} in the affine case:

\begin{theorem}
    \label{maintheorem_aff}
    (1) For any partition $\bla$, $\Phi$ maps the dual PBW basis $|e_{\bla}\rangle^*$ into 
    $K_\sT(\CM_\bd)$ where $\bd = |\bla|$. 

    (2) For any $\bd \in \BZ_{\geq 0}^n$, $\{|e_{\bla}\rangle^*\}_{|\bla|=\bd}$ form a $\BZ_\bu$-basis
    of $K_\sT(\CM_\bd)$.
\end{theorem}

\noindent This theorem will be proved in Section \ref{stab_and_pbw}.

The analog of Corollary \ref{dual_verma_corr} is the following: 

\begin{corollary}
    \label{dual_verma_cor_aff}
    Fix $a_1, \ldots, a_n, l \in \BC$ with $l \neq 0$ and let $q_0 \in \mathbb{C}$ be generic. Let $\BC_a$ denote the evaluation module of $\BZ_\bu$, i.e., the target of the map
    \[ 
        \BZ_\bu = \BZ[u_1^\pm, \ldots, u_n^\pm, q^\pm, p^\pm] \xrightarrow{u_i = q^{a_i}, p = q^l, q = q_0} \BC.
    \]
    Then $K^\text{int} \otimes_{\BZ_\bu} \BC_a$ is isomorphic to the dual Verma module of $U_{q_0}(\widehat{\mathfrak{gl}}_n)$, where $\psi_i$ 
    acts by $q^{i-a_i}$ for $i=1,2,...,n$ and $c$ acts by $q^{n+l}$.
\end{corollary}
 
Corollary \ref{dual_verma_cor_aff} will be proved in Section \ref{action_almost_preserve}.
The main task is to control the denominators appearing in the action, see Proposition \ref{coef_almost_int}. 

\begin{remark}\label{avoid_crit_level}
The condition that $l \neq 0$ is due to the possible denominators in the $f_{[i;j)}$ action.
In terms of the lowest weight, this means that we are away from the critical level where $c$ acts by $q^n$. One cannot specialize to $l=0$ (i.e. $p=1$) because there are denominators of the form $1-p ^{2k}$ in the action of $f_{[i;j)}$. \hide{Note that crit level is $q^n$ not $q^{-n}$. Mainly because we are using lowest weight not highest weight modules.}
\end{remark}

\begin{remark}
Using the RTT integral form of $\glnhat$ as defined in \cite{finkelberg2019shifted}, one can define the integral form of the universal Verma module similar to (\ref{univ_verma}), as well as the integral form of the universal \textit{dual} Verma module. Then the results above say that,
after a mild localization, $K^\text{int}$ can 
be identified with the universal dual Verma module.
\end{remark}

\subsection{Stable envelopes} \label{stab_and_pbw}

In this section, we introduce a variant of the stable envelopes first defined in \cite{maulik2012quantum}, which will be used in the proof of theorem \ref{maintheorem_aff}. To begin with, we collect the ingredients used in the definition of stable envelopes.

\subsubsection{Subtorus \texorpdfstring{$\sA$}{A} and its fixed locus}

Consider the subtorus $\sA = (u_1, \ldots, u_n, p) \subset \sT$ (i.e. $\sA$ is $\sT$ without the $\mathbb{C}^*_q$ factor) 
and fix a 1-parameter subgroup of $\sA$
\begin{equation}\label{1ps_A}
  \sigma_\sA (t) = (t, t^2, \ldots, t^{n}, t^{-n-1})
\end{equation}
We will consider the attracting manifold
with respect to $\sigma_\sA$. 
A weight $p^m u_i^2 / u_j^2$ is attracting if $m<0$ or $m=0$ and $i>j$. It's repelling if $m>0$ or $m=0$ and $i<j$. 

\begin{lemma}
The connected components of the $\sA$-fixed points of $\CM_\bd$ are 
affine spaces, and each contains exactly one point of $(\CM_\bd)^\sT$.
\end{lemma}

\begin{proof}
From the quiver description, the $\sA$-fixed points are given by a product of quivers of the form

\begin{picture}(200,110)(40,-30)\label{fixed_pt_quiver}

    
    
    \put(115,31){\vector(1,0){50}}
    \put(135,34){$Y_{j+1}$}
    
    \put(175,31){\vector(1,0){40}} 
    \put(195,34){$Y_{j+2}$}
    
    \put(247,31){\vector(1,0){38}}
    \put(255,34){$Y_{k}$}
    
    
    \put(110,31){\circle*{10}}
    \put(110,50){\circle{30}}
    \put(95,47){$X_{j}$}
    \put(117,22){$\color{red}{V_{j}}$}
    \put(102,36){\vector(4,-1){5}}
    
    \put(170,31){\circle*{10}}
    \put(170,50){\circle{30}}
    \put(155,47){$X_{j+1}$}
    \put(177,22){$\color{red}{V_{j+1}}$}
    \put(162,36){\vector(4,-1){5}}
    
    \put(225, 30){${\cdots}$}

    \put(290,31){\circle*{10}}
    \put(290,50){\circle{30}}
    \put(275,47){$X_{k}$}
    \put(297,22){$\color{red}{V_{k}}$}
    \put(282,36){\vector(4,-1){5}}
    
    \put(75,-20){\line(1,0){10}}
    \put(75,-10){\line(1,0){10}}
    \put(75,-20){\line(0,1){10}}
    \put(85,-20){\line(0,1){10}}
    \put(88,-20){$\color{blue}{\BC w_{j}}$}
    
    \put(85,-10){\vector(2,3){25}}
    \put(98,4){$A_{i}$}

    \end{picture} 

\noindent 
with the same stability condition and moment map condition as before.
(Note that $k$ could be bigger than $n$. In that case, $V_i$ and $V_{i+n}$ comes from the same $V_i$ in the original quiver, but there is no arrow between them because their weight differ by a non-zero power of $p$.)

Choose a basis $\{v_1,v_2,...,v_{\dim V_j}\}$ of $V_j$ so that $X_j$ has the form
\[
    \begin{pmatrix}
    0 & & & & & *\\
    1 & 0 & & & & *\\
     & 1 & & & & *\\
     & & \ddots & & & \vdots\\
     & & & \ddots & & \vdots\\
     & & & & 1 & *\\
    \end{pmatrix}
\]
For each $i>j$, take the vectors 
\[ 
    Y_i Y_{i-1}... Y_{j+1}v_t, \, t=1,2,...,\dim V_i
\]
to be a basis of $V_i$. (The stability condition and the moment map condition ensures that this is a basis of $V_i$.) And all other arrows in this basis are determined by $X_j$. So the only free parameters are the entries marked by $*$ in the above matrix, which form an affine space.

\end{proof}

We will denote by $\BA_\bla$ the connected component of $\left(\CM_\bd \right)^{\sA}$ that contains $\bla$. In the general definition of stable envelopes, one frequently consider the restriction of  K-theory classes to the $\sA$-fixed locus. In our case, this is determined by the restriction to the $\sT$-fixed points. So we will restrict to the $\sT$-fixed points instead to make the notations less cumbersome.

\subsubsection{Semi-polarization}
Defining stable envelopes usually requires a choice of polarization, which is roughly half of the tangent bundle. $\CM_\bd $ does not have a polarization, but we will see that an analog of stable envelopes can still be defined. The K-theory class $\rV$ below will play the role of polarization.

Let
\[ 
    \rV = \bigoplus_{i=1}^n \CA_i, \,\,
    \mathscr{U} = \bigoplus_{i=1}^n \CB_i 
\] 
where $\CA_i, \CB_i$ denote the tautological bundles from $A_i, B_i$ in the quiver. 


\begin{lemma} \label{v_equals_attr}
For each fixed point $\bla \in \CM_\bd$, viewing $\rV, \mathscr{U}, T \CM_\bd $ as elements in  $K_\sA(\CM_\bd)$, we have 
\begin{align*}
    \left(\rV | _\bla \right)_{{mov}} &= (T\CM_\bd|_\bla)_{rep}, \\
    \mathscr{U}|_\bla &= (T\CM_\bd|_\bla)_{attr},
\end{align*}
where the subscript $mov$ means the moving (non-fixed) part {with respect to the torus $\sA$}, and $rep/attr$ means repelling/attracting directions.
\end{lemma}

\begin{proof}
When only considering the $\sA$-action (setting $q=1$), the formula (\ref{tangentbundle}) simplifies to
\[ 
    \sum_{i=1}^n\frac{\CV_i}{u_i^2}+ \sum_{i=1}^n\frac{u_{i+1}^2}{\CV_i}
\] 
which equals $\rV + \mathscr{U}$ termwise. 
$\mathcal{V}_i|_\bla$ is expressed in \eqref{eqn:V_sum_box}. 
For a color $i$ box $\sq$, if it is in $\lambda ^k$ with $k \leq i$, then $\chi _\sq$ is $u_k^2$ times a non-negative power of $p$, so $\chi_\sq / u_i^2$ is either repelling or fixed; 
if $k > i$, then $\chi_\sq$ is $u_i^2$ times a \textit{positive} power of $p$, so $\chi_\sq / u_i^2$ is repelling. This proves the first equality. The second equality can be proved similarly -- if a color $i$ box is in $\lambda ^k$ with $k \geq  i+1$ then $\chi _\sq$ must contain positive power of $p$, 
which makes $u_{i+1}/\chi _\sq$ attracting.
\hide{Also note that for $i=n$, the term is $u_1^2p^{-2}/\CV_n$ which always have negative power of $p$.}
\end{proof}

\begin{remark}
Note that the above lemma is not true if we consider these classes as $\sT$-equivariant. 
This is analogous to the fact that when considering stable envelopes for Nakajima quiver varieties (see e.g. \cite{maulik2012quantum}), one considers the subtorus $\sA \subset \sT$ where $\sT / \sA \simeq \mathbb{C}^*$ is the torus that scales the symplectic form. 
\end{remark}


\subsubsection{Attracting Set and order of fixed points}

We write $\bmu \succ \bla$ if $\bmu$ is in the full attracting set of $\BA_\bla$. (See \cite{maulik2012quantum} for the exact definition of full attracting set. Roughly speaking,
this means taking the attracting locus and taking closure iteratively.) This defines a partial ordering on the fixed points. 
We also choose a total ordering on the fixed points such that $\bmu > \bla$ if $h(\bmu) > h(\bla)$, where $h$ is the sum of the height (row index) of boxes in $\bmu$ or $\bla$. If $h(\bla) = h(\bmu)$ then they can be in either order.
For example, the partition with only one row in each partition is the smallest for either $>$ or $\succ$. It's attracting set is an open subset of $\CM_\bd $.

The following lemma will be used in studying stable envelopes, c.f. Lemma \ref{lem:stab_bounded}.

\begin{lemma}\label{lem:color_compare}
Let $h_c(\bla)$ denote the sum of height of color $c$ boxes in a partition $\bla$.
If $\bmu \succ \bla$, then for every color $c \in \{1,2,...,n\}$, $h_c(\bmu) \geq  h_c(\bla)$, and equality holds only when the color $c$ boxes in $\bmu$ and $\bla$ are in the same positions. In particular, $\bmu \succ \bla$ implies $h(\bmu) > h(\bla)$.
\end{lemma}

\begin{proof}
This can be proved in the same way as Proposition \ref{prop:attr_vanish}. The only difference is that, instead of comparing the number of color $c$ boxes in a given partition $\mu ^k$ or $\lambda ^k$, we need to compare the number of color $c$ boxes in a given row $l$. 
Note that a color $c$ box in row $l$ in $\bmu$ necessarily lies in the partition $\mu ^{\overline{c-l}}$. 
If two color $c$ boxes $\sq_1, \sq_2$ are in row $k_1$ and row $k_2$ respectively, and $k_1 > k_2$, then the weight $\chi _{\sq_1}/ \chi_{\sq_2}$ is repelling. 
So the same argument analyzing the weight in local coordinates $U_\bla$ as in Proposition \ref{prop:attr_vanish} also applies here.
\end{proof}

\subsubsection{Degree in \texorpdfstring{$\sA$}{A}}

\begin{definition}
For an element $\alpha \in K_\sT(pt)$, 
viewing $\alpha$ as a Laurent polynomial in variables $u_1,...,u_n, p$ with coefficients in $\mathbb{Z}[q, q ^{-1}]$, the $\sA$-degree of $\alpha$, denoted by $\deg_\sA \alpha$, is defined as 
the Newton polygon of $\alpha$ in the character lattice of $\sA$. 
In other words, let 
\[ 
    S = \{w \in \text{char}(\sA) \mid \text{coefficient of $w$ in $\alpha$ is not 0}\},
\]
then $\deg_\sA \alpha$ is the convex hull of $S$ in $\text{char}(\sA) \otimes \BR$.
\end{definition}

\subsubsection{Stable envelopes}
Let
\[ 
    L :=\left( \prod_{i=1}^n \det \CV_i \right)^\epsilon
\]
where $\epsilon > 0$
is a very small real number. (We will explain what power $\epsilon$ means. $\epsilon $ needs to be very close to zero and may depend on $\bd $. The exact bound will be clear from the proof.)

\begin{proposition}
[\cite{okounkov2021inductive}, Theorem 2]
\label{stab_exists}
For each $\bla$, there exists $s_\bla \in K_\sT(\CM_\bd)$ such that:
\begin{itemize}
  \item[(i)] $s_\bla$ is supported on the full attracting set of $\BA_\bla$.
  \item[(ii)] $s_\bla|_\bla = \Lambda^\bullet\left(\rV_\bla^{mov}\right)^\vee$.
  \item[(iii)] For any $\bmu \succ \bla$, 
  \[ 
    \deg_\sA s_\bla|_\bmu + \mathrm{weight}_\sA L|_\bla \subset \deg_\sA \Lambda^\bullet \left(\rV_\bmu^{mov}\right)^\vee + \mathrm{weight}_\sA L|_\bmu.
  \]
\end{itemize}
\end{proposition}

For a line bundle $L'$, $\mathrm{weight}_\sA L'|_\bla$ is a vector in $\text{char}(\sA) \otimes \BR$. The definition of $\mathrm{weight}_\sA$ also applies to an $\epsilon $ power of a line bundle by multiplying this vector by $\epsilon$. So it applies to the $L$ defined above. And $+ \mathrm{weight}_\sA L|_\bla$ means translating the polygon by this vector.

\begin{remark}
{\cite{okounkov2021inductive} proved this under the simplifying assumption that there exists a polarization $T ^{1/2}$. But it can be seen that $T ^{1/2}$ can be replaced by $\rV$ in our setting and the proof there goes through -- the property that is used in the proof is that for each fixed point $\bla$, the Newton polygon $\deg_\sA \Lambda ^\bullet \left(\rV_\bla^{mov}\right) $ agrees with $ \deg_\sA \Lambda ^\bullet \left(T \CM_\bd |_\bla \right)_{rep}$ up to translation. In fact, relaxing the condition of polarization to something that resembles the repelling direction is one of the main goals of \cite{okounkov2021inductive}.

Also note that \cite{okounkov2021inductive} used topological K-theory instead of algebraic K-theory. But since $\CM_\bd $ has an affine cell decomposition (\cite{neguct2018affine} Remark 3.10), there is no difference between topological and algebraic K-theory.
\hide{I think Andrei's argument relies on that every fixed point, the repelling direction should link to another fixed point (or the fp is maximal). Not sure how this is proved even for Nakajima quiver variety.}
}
\end{remark}




\subsection{Proof of the main theorem}

The same argument as in Section \ref{proof_finite_case} shows that part (1) of the theorem implies part (2). To establish part (1), 
it suffices to find K-theory classes whose pairing with the PBW basis is an upper triangular matrix with invertible diagonal entries.
This is where stable envelopes come in handy.

We will leverage the rigidity argument introduced in \cite{okounkov2015lectures} to compute 
the inner product of stable envelopes and PBW basis. 
The logic can be summarized as follows: 
first, by the localization formula (\ref{pairing_on_fp}), the inner product of K-theory classes can be computed by summing over fixed points. The degree constraint (property (iii) in the definition of $s_\bla$, and Lemma \ref{pbw_bounded} below) implies that the contribution from each fixed point is ``bounded'' in the sense defined below. 
Hence, the left-hand side of (\ref{pairing_on_fp}) also has this property. The properness of correspondences implies that the LHS is a Laurent polynomial. 
Together with boundedness, this implies that the LHS is constant in $\sA$. Thus, it can be computed by sending equivariant variables in $\sA$ to any limit, for which we can choose a convenient one to make most fixed-point contributions vanish; see Proposition \ref{compute_pairing}.

For a rational function $f \in \text{Frac}\,\BZ_\bu$ and a 1-parameter subgroup $\xi (t): \mathbb{C}^* \to \sA$, define
\[ 
    \lim _{\xi (t)\to \infty }f := \lim _{t \to \infty }f| _{(u_1,...,u_n, p) = \xi (t)}.
\]
(Note that $q$ is not replaced in the above definition. So $\lim _{\xi (t)\to \infty }f$ is a function of $q$ if the limit exists.)
We say that $f$ remains bounded in any limit of $\sA$ (or bounded for short) if $\lim _{\xi (t)\to \infty }f$ exists for any subgroup $\xi (t)$. 

Recall the subgroup $\sigma_\sA (t)$ defined in formula (\ref{1ps_A}). 

\begin{lemma}\label{lem:stab_bounded}
Given fixed points $\bla, \bmu$, let
\[ 
    f _{\bla \bmu} = \frac{s_\bla|_\bmu}{\Lambda ^\bullet (T ^{<0}_\bmu)^\vee }.
\]
Then:
\begin{itemize}
\item[(1)] If $\bmu < \bla$, then $f _{\bla \bmu} = s_\bla|_\bmu = 0$.
\item[(2)] If $\bmu > \bla$, then $f _{\bla \bmu}$ is bounded in any limit of $\sA$,
and $\lim _{\sigma _\sA(t) \to \infty }f _{\bla \bmu}=0$. 
\item[(3)] If $\bmu=\bla$, then $f _{\bla \bmu}=1$.
\end{itemize}
\end{lemma}
\begin{proof}
Statements (1) and (3) follows directly from the definition of $s_\bla$. 
The boundedness in property (2) follows from condition (iii) of $s_\bla$.

Furthermore, note that 
condition (iii) of $s_\bla$ combined with Lemma \ref{lem:color_compare} shows that  
\[ 
    \lim _{\sigma _\sA(t) \to \infty }f _{\bla \bmu} \cdot t ^{k \epsilon }
\]
exists for some positive integer $k$. Here the factor $t ^{k \epsilon }$ comes from the shift $\mathrm{weight}_\sA L|_\bla  - \mathrm{weight}_\sA L|_\bmu$. So $\lim _{\sigma _\sA(t) \to \infty }f _{\bla \bmu}=0$

\hide{I was wondering whether there's a more generic argument to say $f_{\bla \bmu}->0$ without using too much information about $L$ and $\sigma _\sA$. For stable envelope for Nakajima quiver variety this is ok because one can choose a direction such that $f->0$. But here we need a specific $\sigma _\sA$ limit to match the limit in $|e_{\bla}\rangle$. So a careful control of $L$ is unavoidable. Another thought: maybe this can be done by arguing that $L$ is an ample line bundle? But probably that's the same as our proof here. (Also, one needs to show that a fixed point in the attracting set can be connected to by a chain of P1's. I showed this in an appendix before.)}
\end{proof}

Now, consider the K-theory classes
\begin{equation}
    \label{pbw_normalize}
    \widetilde{|e_{\bla}\rangle} := \left. |e_{\bla}\rangle \middle/  \prod_{i=1}^n u _{i+1}^{d_i} \right.
\end{equation}
where $u _{n+1} := u_1 p ^{-2}$.
\begin{lemma}\label{pbw_bounded}
Given fixed points $\bla, \bmu$, let
\[ 
    g _{\bla \bmu} = \frac{\widetilde{|e_{\bla}\rangle}|_\bmu}{\Lambda ^\bullet (T ^{> 0}_{\bmu})^\vee }.
\]
Then $g _{\bla \bmu}$ remains bounded in any limit of $\sA$. Furthermore, if $\bmu < \bla$, then $\lim _{\sigma _\sA(t)\to \infty } g _{\bla \bmu}=0$.
\end{lemma}

\begin{proof}
Recall that $|e_{\bla}\rangle$ is constructed by applying $e_{[i;j)}$
coming from the columns of $\bla$. More explicitly, suppose that $\bla$ has $k$ columns and
\[ 
    e_\bla = \prod_{l=1}^k e_{[i_l;j_l)},
\]
Then
\[ 
    \langle \bmu | e_\bla \rangle  = \sum_{\bmu_1 = \bmu, \bmu_2, ..., \bmu_k, \bmu_{k+1}=\emptyset}
    ^{\substack{\bmu_{l+1} \subset \bmu_l \text{ and }\\
     \bmu_l \backslash \bmu_{l+1} \text{ has boxes of color }[i_l;j_l)}} 
    \prod_{l=1}^k 
    \langle \bmu_l | e_{[i_l;j_l)} | \bmu_{l+1} \rangle .
\]
Recall that 
\begin{equation}\label{eqn:gblabmu}
    g _{\bla \bmu} = \frac{\langle \bmu| e_\bla \rangle }{\Lambda ^\bullet (T ^{> 0}_{\bmu})^\vee  \cdot \prod_{i=1}^n u _{i+1}^{d_i}},
\end{equation} 
We will show that each $\langle \bmu_l | e_{[i_l;j_l)} | \bmu_{l+1} \rangle$ can be paired with suitable terms from the denominator in \eqref{eqn:gblabmu} to make it bounded in $\sA$.

From formula \eqref{actioncoeff}, each $\langle \bmu_l | e_{[i_l;j_l)} | \bmu_{l+1} \rangle$ is a product of $R_{ij}^+$ and a $\zeta\tau_+$ term for each box in $\bmu_{l}\backslash \bmu_{l+1}$.

We first analyze the $\zeta \tau _+$ terms. 
Let $\bnu = \bmu_{l+1}$.
Formula (\ref{simple_zeta_tau}) can be rewritten as
\begin{multline}\label{split_ground}
    \zeta\left(\frac{z}{\chi _\bnu}\right)\tau_+(z) =\\
    \underbrace{\frac{
        \prod_{\sq \in \bnu \text{ addable, not on the bottom row}}^{\col \sq = \col z + 1} 
        \left(\frac{\sqrt{\chi _\sq}}{q}- \frac{zq}{\sqrt{\chi _\sq}}\right)
    }{
        \prod_{\sq \in \bnu \text{ removable}}^{\col \sq = \col z}\left(\sqrt{\chi _\sq} - \frac{z}{\sqrt{\chi _\sq}}\right)
    }}_{\text{(I)}}
    \cdot 
    \underbrace{\left(1-\frac{z q^2}{\chi _\boxtimes }\right)
    \vphantom{\frac{\prod \frac{1}{\sqrt{x_\sq}}}{\left(\prod \frac{1}{\sqrt{x_\sq}}\right)}} 
    }_{\text{(II)}}
    \cdot 
    \underbrace{\frac{\sqrt{\chi_\boxtimes}}{q} 
    \vphantom{\frac{\prod \frac{1}{\sqrt{x_\sq}}}{\left(\prod \frac{1}{\sqrt{x_\sq}}\right)}} 
    }_{\text{(III)}}
\end{multline}
where $\boxtimes$ stands for the addable box on the bottom row of the ($i+1$)-th partition in $\bnu$ where $i = \col z$. In particular, $\chi _\boxtimes = u _{i+1}^2 q ^{2\bnu^{i+1}_0}$. Thus, part (III) in (\ref{split_ground}), when taking $z$ to be each $\bsq \in \bmu_l \backslash \bmu_{l+1}$ and letting $l$ ranging from 1 to $k$, 
cancels the product 
\[ 
    \prod_{i=1}^n u _{i+1}^{d_i}
\]
up to powers of $q$.

Part (I) in (\ref{split_ground}) is obviously bounded since the removable box in color $i$ and the addable box in color $i+1$ come in pairs (unless the color $i+1$ box is on the bottom row, and that's why we separate it). 

For part (II), from the proof of Lemma \ref{v_equals_attr}, we know that for each $\bsq \in \bmu$ with $\col \bsq = i$, there is a corresponding attracting direction with weight
\[ 
    w = \frac{u _{i+1}^2}{\chi _{\bsq }q^{m} }
\]
for some number $m$. Then 
\[ 
    \frac{1-\frac{\chi _\bsq q^2}{\chi _\boxtimes }}{1-w ^{-1}}
\]
is bounded. So part (II), when taking $z$ to be each $\bsq \in \bmu_l \backslash \bmu_{l+1}$ and letting $l$ ranging from 1 to $k$,
then dividing by $\Lambda ^\bullet (T ^{> 0}_{\bmu})^\vee$, is bounded in any limit of $\sA$. 

Now we examine the $R_{ij}^+$ factor.
In formula \eqref{rij}, both $1/(1-z _{i+1}/z_i q^2)$ and the $\zeta $ factors are obviously bounded. So we proved that $g_{\bla \bmu}$ is bounded in any limit of $\sA$.


When $\bmu < \bla$, there exists $l$ such that when going from $\bmu_{l+1}$ to $\bmu_l$, at least one added box is put at a lower row than its standard position. (Standard position means that the boxes of colors $[i_l;j_l)$ are added to the $i_l$-th partition in row $0,1,..., j_l- i_l -1$ respectively.) 
Suppose the first such box corresponds to number $m$ in $[i_l;j_l)$, then $z _{m}/z_{m-1} $ in $R_{ij}^+(\bmu_l \backslash \bmu_{l+1})$ goes to $\infty $ in the limit ${\sigma _\sA(t)\to \infty }$. 
\hide{The symmetrization won't affect this argument. A better way to say is that, if $z_m, z_{m-1}$ are both in standard position, then their ratio is a power of $q$. Moving boxes lower makes it go to $\infty $ faster, and higher slower. We can always find $m$ such that $z_m$ is lower and $z_{m-1}$ is higher or unchanged.}
This makes each product in $g _{\bla \bmu}$ (and hence $g_{\bla \bmu}$ itself) go to zero in the limit. 
\end{proof}

Recall the pairing on $K_\sT(\CM_\bd )$ defined in (\ref{shapo_k}). Let $\widetilde{s}_\bla = s_\bla \otimes \mathcal{D}_\bd $ to compensate for the twist in the pairing.
\begin{lemma}
For any fixed point $\bla, \bmu$, the pairing $(\widetilde{s}_\bla, \widetilde{|e_{\bmu}\rangle}) \in \BZ[q, q ^{-1}]$, i.e., it is a Laurent polynomial that only depends on $q$.
\end{lemma}
\begin{proof}
By lemma \ref{corr_proper}, $(\widetilde{s}_\bla, \widetilde{|e_{\bmu}\rangle}) \in \BZ_\bu$. By the localization formula,
\begin{equation}\label{pairing_on_fp}
(\widetilde{s}_\bla, \widetilde{|e_{\bmu}\rangle}) = \sum_{\text{fixed point }\bnu} f _{\bla \bnu}g _{\bmu \bnu}.
\end{equation}
By lemmas \ref{lem:stab_bounded} and \ref{pbw_bounded}, each term on the right-hand side is bounded in any limit of $\sA$, thus so is the left-hand side. But the LHS is a Laurent polynomial in $\BZ_\bu$, thus it must not depend on $u_1,...,u_n, p$, hence it is an element in $\BZ[q , q ^{-1}]$.
\end{proof}

\begin{proposition} \label{compute_pairing}
For $\bmu > \bla$, $(\widetilde{s}_\bla, \widetilde{|e_{\bmu}\rangle})=0$. For $\bmu = \bla$, $(\widetilde{s}_\bla, \widetilde{|e_{\bla}\rangle})=q^N$ for some integer $N$.
\end{proposition}
\begin{proof}
Since $(\widetilde{s}_\bla, \widetilde{|e_{\bmu}\rangle})$ does not depend on $u_1,...,u_n, p$, it can be computed in any limit of $\sA$.

For $\bmu > \bla$, for any fixed point $\bnu$, either $\bmu > \bnu$ or $\bnu > \bla$ (or both), hence by lemmas \ref{lem:stab_bounded} and \ref{pbw_bounded}, either $\lim _{\sigma _\sA(t)\to \infty } g _{\bmu \bnu} = 0$ or $f _{\bla \bnu}=0$. So the RHS of (\ref{pairing_on_fp}) is 0.

For $\bmu = \bla$, 
we compute $(\widetilde{s}_\bla, \widetilde{|e_{\bla}\rangle})$ in the limit $\sigma _\sA(t)\to \infty $. The same argument shows that the only contribution comes from $f _{\bla \bla}g _{\bla \bla}$ after taking the limit. 
Also, the symmetrization at the beginning of $R^+$ in formula (\ref{rij}) can be dropped because that will make boxes in non-standard positions. (See the last part of the proof of Lemma \ref{pbw_bounded}.) 

By construction, $f_{\bla \bla}=1$. So we only need to consider $g_{\bla \bla}$ in the limit. 
For a rational function of the form
\[ 
    \frac{1- \frac{\chi _\sq}{\chi _\bsq}q^k}{1-\frac{\chi _\sq}{\chi _\bsq}}
\]
for some integer $k$. If $\sq, \bsq$ are not in the same partition $\la ^i$ of $\bla$, then the limit is either 1 or $q^k$, and can be dropped for the purpose of this proof. The non-trivial factors come from the following. (Again, we are focusing on $e_{[i;j)}$ coming from one column of $\bla$ as in the proof of lemma \ref{pbw_bounded}, and note that all boxes must land in standard position in the $i$-th partition, and when applying it there is no box above the $(j-i)$-th row): 
\begin{itemize}
\item The term in (\ref{split_ground}) for $z=\chi _\bsq$ is only non-trivial when the removable box is next to $\bsq$ and the addable box is one row above $\bsq$, this results in a factor of 
\[  
    \frac{1- q ^{2(\lambda ^i_m - \lambda ^i _{m+1}+1)}}{1-q^2}
\]
where $m = \col \bsq - i$.
\item In formula (\ref{actioncoeff}), each $\bsq$ comes with a factor of $q-q ^{-1}$, this cancels the denominator above, up to powers of $q$.
\item In $R^+$, for adjacent boxes $\sq, \bsq$ with $\col \sq = m-i+1, \col \bsq = m-i$,
\[ 
    \frac{1}{1- \frac{\chi _\sq}{\chi _\bsq q^2}} = \frac{1}{1-q ^{2(\lambda ^i _{m+1}-\lambda ^i_m - 1)}}
\]
This almost cancels the numerator in the first item, except for 
one term from $m=j$. (Note that $\lambda ^i _{j+1}=0$ because we 
apply columns from right to left.) 
\item The pure $q$ weight in $\Lambda ^\bullet (T_\bla ^{\geq 0} )^\vee $:
It can be computed from formula (\ref{tangentbundle})  that each $e_{[i;j)}$ will increase this by $1- q ^{2( \lambda ^i _j +1)}$, this cancels the remaining term from the first item. 
\end{itemize}
So we see that $(\widetilde{s}_\bla, \widetilde{|e_{\bla}\rangle})$ is indeed a power of $q$.
\end{proof}


Therefore, the pairing between $\widetilde{s}_\bla$ and $\widetilde{|e_{\bmu}\rangle}$,  
when $\bla$ and $\bmu$ range over all fixed points, forms a lower triangular matrix  
with powers of $q$ on the diagonal. Consequently, its inverse lies in $\BZ_\bu$.  

This implies that the dual PBW basis  are $\BZ_\bu$-linear combinations  
of the $\widetilde{s}_\bla$'s. Hence, these dual PBW basis belong to $K_\sT(\CM_\bd)$.  
This completes the proof of Theorem \ref{maintheorem_aff}.

\subsection{\texorpdfstring{$\glnhat$ action almost preserves $K_\sT(\CM_\bd)$}{glnhat action almost preserves K}}
\label{action_almost_preserve}

The action of $e_{[i;j)}$ and $\tau(e_{[i;j)})$ preserves $K^\text{int}$,  
see Remark \ref{rmk:tau_proper}.  
However, the operator $f_{[i;j)}$ does not always preserve $K^\text{int}$.  
(For example, this happens for $n=2$ and the operator $f_{[1;3)}$ from $K_\sT(\CM_{(1,1)} )$ to $K_\sT(\CM_{(0,0)})$. The correspondence $\ofZ_{[1;3)}$ is isomorphic to the blow up of $\mathbb{P}^1 \times \mathbb{C}$ at one point.)
Nonetheless, we can control the denominators in the $f_{[i;j)}$ action,  
which leads to a proof of Corollary \ref{dual_verma_cor_aff}.  

For a given $[i;j)$, consider the element $\tau(f_{[i;j)})$ and write it as a $\BQ(q, p)$-linear combination  
of PBW basis $e_\bla$, i.e.,  

\begin{equation}
    \label{def_c_bla}
    \tau(f_{[i;j)}) = \sum_{|\bla|=[i;j)} c_\bla e_\bla,
\end{equation}  
\noindent where the coefficients $c_\bla$ are, a priori, elements of $\BQ(q, p)$.  
We will control the denominators of $c_\bla$ and this in turn controls the denominator of the $f_{[i;j)}$ action (because $\tau (e_\bla)$ preserves $K ^{^\text{int}}$).  

Let  
\begin{multline}\label{loc_coef_def}
    \BZ[q^\pm, p^\pm]_{1-p^\#, 1-q^\#} :=  
    \left\{\ \frac{f}{g}\
    \bigg|\
    \substack{f \in \BZ[q^\pm, p^\pm], \, 
    g = \prod_{i=1}^m (1 - p^{2k_i})\prod_{j=1}^h (1 - q^{2l_j}) \\
    \text{ for some integers } k_1, \dots, k_m, l_1,...,l_h} \right\}.
\end{multline}

\begin{proposition} \label{coef_almost_int}
The coefficients $c_\bla$ live in $\BZ[q^\pm, p^\pm]_{1-p^\#, 1-q^\#}$.
\end{proposition}  

This will be proved in section \ref{prove_denom}.  
This implies that $f_{[i;j)}$ can be written as a sum of products of $\tau(e_{[i';j')})$ whose coefficients have denominators no worse than a product terms of the form $1-p^{2k}$ and $1-p ^{2k'}$. Hence, the action of $\glnhat$ preserves $K^\text{int}_{1-p^\#, 1-q^\#}$, where the subscript means inverting elements of the form $1-p^{2k}, 1-q ^{2k'}, k, k' \in \BZ$. This combined with theorem \ref{maintheorem_aff} proves corollary \ref{dual_verma_cor_aff}.

\subsubsection{Bilinear pairing}
 
There is a non-degenerate bilinear pairing between $\glnhatleq$ and $\glnhatgeq$ defined by  
\begin{align}
    \langle f_i, e_j \rangle  &= (q - q^{-1})\delta_j^i, \label{eqn:fe_pair} \\
    \langle P_k, P_l \rangle &= \delta_{k+l}k \cdot \frac{  
        (q^{nk} - q^{-nk})  
    }{  
        (p^k - p^{-k})(q^{nk}p^k - q^{-nk}p^{-k}),  
    }  \label{eqn:pp_pair}
\end{align} 
on generators and satisfies  
\[ 
    \langle a  a', b \rangle = \langle a \otimes a', \Delta(b) \rangle,
\]  
where $\Delta $ is the coproduct. (We omit the pairing on the Cartan part. They are all constants not depending on $q,p$.)

From these, we can deduce that the pairing  
\begin{equation}
    \label{s+s-pairing}
    \langle  f_{i'_1}\cdots f_{i'_{m'}}P_{-j'_1}\cdots P_{-j'_{l'}} , \; e_{i_1}\cdots e_{i_m}P_{j_1}\cdots P_{j_l}\rangle 
\end{equation}  
(assuming that the $j$'s and $j'$'s are in increasing order since the $P$'s commute) is non-zero only when $m = m'$, $l = l'$, $j_i = j'_i$ for $i=1,\dots,l$, in which case the pairing is equal to
\begin{equation}\label{eqn:ss-pairing-res}
    \prod_{i=1}^l \langle P_{-j_i}, P_{j_i} \rangle \cdot \text{a rational function in $q$}.
\end{equation}

Pairing between root generators and the $P_k's$ can be deduced from the embedding of $\glnhat$ into the shuffle algebra (Proposition 3.33 of \cite{neguct2013quantum}) and invoking results about shuffle algebra.
\cite{neguct2013quantum}, Proposition 4.2 implies that  
\begin{equation}
    \langle f_{[i;i+nk)}, P_k \rangle = (-q^{-1})^{nk} (-1)^{k-1} \frac{q^k-q^{-k}}{p^k-p^{-k}}
\end{equation}  
\begin{equation}
    \langle P_{-k}, e_{[i;i+nk)} \rangle = (-1)^{k-1} \frac{q^k-q^{-k}}{p^k q^{nk} - p^{-k} q^{-nk}}
\end{equation}  
Combined with Lemma 3.20 and Proposition 3.25 of \loccitt, which simplifies the pairing between a root generator and a product of several elements, 
we see that  
\begin{align}
    \langle f_{[i;j)}, \; e_{i_1}\cdots e_{i_m}P_{j_1}\cdots P_{j_l} \rangle &= \prod_{i=1}^l \frac{1}{p^k - p^{-k}} \cdot \text{rational function of $q$}, \label{f_pair_with_plus} \\  
    \langle f_{i_1}\cdots f_{i_m}P_{-j_1}\cdots P_{-j_l}, \; e_{[i;j)} \rangle &= \prod_{i=1}^l \frac{1}{p^k q^{nk} - p^{-k} q^{-nk}} \cdot \text{rational function of $q$}. \label{e_pair_with_minus} 
\end{align}

\subsubsection{Proof of Proposition \ref{coef_almost_int}}\label{prove_denom}
To prove proposition \ref{coef_almost_int}, we first prove a weaker result:
\begin{proposition}
Let $c_\bla$ be as defined in (\ref{def_c_bla}). Then $c_\bla$ lives in $\BQ(q)[p^\pm]_{1-p^{\#}}$, where the subscript $1-p^{\#}$ is defined in the same way as (\ref{loc_coef_def}).
\end{proposition}

\begin{proof}
Let $S_+$ (resp. $S_-$) be a basis of $\uup$ (resp. $\uum$) whose elements are all of the form $e_{i_1}\cdots e_{i_m}P_{j_1}\cdots P_{j_l}$ (resp. $f_{i_1}\cdots f_{i_m}P_{-j_1}\cdots P_{-j_l}$).

Write 
\[ 
    f_{[i;j)} = \sum_{s \in S_-} b_s \cdot s.
\]
We claim that $b_s \in \BQ(q)[p^\pm]$ for all $s \in S_-$. Indeed, $b_s$ can be computed from the pairing $\langle f_{[i;j)}, s \rangle$ for all $s' \in S_+$ as in (\ref{f_pair_with_plus}).

The pairing between $S_-$ and $S_+$ can be arranged into a block diagonal matrix 
\[
\begin{pmatrix}
{M_{J_1}} &  &  &  \\
 & {M_{J_2}} &  &  \\
&& \ddots & \\
 &  &  & {M_{J_m}}
\end{pmatrix}
\]
where for a given index set $J_i = {j_1,...,j_{k_i}}$, 
$M_{J_i}$ is a square matrix corresponding to the pairing between elements in $S_-$ whose ``$P$ part'' is $P_{-j_1}... P_{-j_{k_i}}$ and elements in $S_+$ whose ``$P$ part'' is $P_{j_1}... P_{j_{k_i}}$.
Formula \eqref{eqn:ss-pairing-res} implies that each $M_{J_i}$ can be written as a matrix that only depends on $q$ times a constant. (This constant comes from pairing between $P_{-j_1}... P_{-j_{k_i}}$ and $P_{j_1}... P_{j_{k_i}}$, which can be computed using \eqref{eqn:pp_pair}.)

Putting things together, and note that denominators in \eqref{eqn:pp_pair} kills the denominators in \eqref{f_pair_with_plus}, we deduce that for  
$s = f_{i_1}\cdots f_{i_m}P_{-j_1}\cdots P_{-j_l}$, 
\[ 
    b_s \in \prod_{i=1}^l (p^{j_i} q^{n j_i} - p^{-j_i} q^{-n j_i}) \cdot \BQ(q),
\]
and in particular, is an element of $\BQ(q)[p^\pm]$.

Our next task is to write $P_k$ in terms of $e_{[i;j)}$. Similar to the above, write
\[ 
    e_{[i;i+nk)} = \sum_{s \in S_+} a_s \cdot s.
\]
(Note that $S_+$ must contain $P_k$ for dimension reason.) The same argument as above shows that for 
$s = e_{i_1}\cdots e_{i_m}P_{j_1}\cdots P_{j_l}$, 
\[ 
    a_s \in \prod_{i=1}^l (p^{j_i} - p^{-j_i}) \cdot \BQ(q).
\]
And the coefficient for $s=P_k$ is $p^k - p^{-k}$.
This shows (by induction on $P_i$'s) that $P_k$ can be written as a sum of products of $e_{[i;j)}$'s where the coefficients are in $\BQ(q)[p^\pm]_{1-p^{\#}}$.

Combining all the above, and noting that 
\[ 
    \tau(f_i) = e_i, \quad \tau(P_{-k}) = P_k,
\]
we deduce that $\tau(f_{[i;j)})$ can be written as a sum of products of $e_{[i;j)}$'s where the coefficients are in $\BQ(q)[p^\pm]_{1-p^{\#}}$.
\hide{Since $G_k$ is scaled by powers of $p$, strictly speaking there is another step that turn $P_k$ to $G_k$, scale, then turn back into $P_k$. But the denominator should be the same.}

Since any product of $e_{[i;j)}$'s can be written as a $\BQ(q)$-linear combination of the PBW basis $e_\bla$, this finishes the proof of the proposition. 
\end{proof}

\begin{proof}[Proof of proposition \ref{coef_almost_int}]
Note that 
\[ 
    c_\bla = (\tau(f_{[i;j)}) \cdot |0 \rangle , \; |e_\bla\rangle^*) = (|0 \rangle , f_{[i;j)}  |e_\bla\rangle^*),
\]
where $|0 \rangle = [\CO_{\CM_0}]$ is the structure sheaf in degree $\bd=0$. So we only need to prove that $f_{[i;j)}  |e_\bla\rangle^*$
has the expected denominator. From the localization formula, all factors of the denominator of $f_{[i;j)}  |e_\bla\rangle^*$ must have the form $1-w$ for some monomial $w$. On the other hand, we have shown that it's an element of $\BQ(q)[p^\pm]_{1-p^{\#}}$. So the factors of the denominator can only be of the form $1-p ^{2k}, 1- q ^{2l}$.

\end{proof}


\appendix

    \bigskip
    \footnotesize
  
    \textsc{Department of Mathematics, Columbia University,
      New York, NY, USA}\par\nopagebreak
    \textit{E-mail address}: \texttt{cs3627@columbia.edu}


\begin{bibdiv}
    \begin{biblist}
    
    \bib{bialynicki1973some}{article}{
          author={Bialynicki-Birula, Andrzej},
           title={Some theorems on actions of algebraic groups},
            date={1973},
         journal={Annals of mathematics},
          volume={98},
          number={3},
           pages={480\ndash 497},
    }
    
    \bib{braverman2005finite}{article}{
          author={Braverman, Alexander},
          author={Finkelberg, Michael},
           title={Finite difference quantum toda lattice via equivariant k-theory},
            date={2005},
         journal={Transformation Groups},
          volume={10},
           pages={363\ndash 386},
    }
    
    \bib{bullimore2016vortices}{article}{
          author={Bullimore, Mathew},
          author={Dimofte, Tudor},
          author={Gaiotto, Davide},
          author={Hilburn, Justin},
          author={Kim, Hee-Cheol},
           title={Vortices and vermas},
            date={2016},
         journal={arXiv preprint arXiv:1609.04406},
    }
    
    \bib{chriss1997representation}{book}{
          author={Chriss, Neil},
          author={Ginzburg, Victor},
           title={Representation theory and complex geometry},
       publisher={Springer},
            date={1997},
          volume={42},
    }
    
    \bib{feigin2011gelfand}{article}{
          author={Feigin, Boris},
          author={Finkelberg, Michael},
          author={Frenkel, Igor},
          author={Rybnikov, Leonid},
           title={Gelfand--tsetlin algebras and cohomology rings of laumon spaces},
            date={2011},
         journal={Selecta Mathematica},
          volume={17},
           pages={337\ndash 361},
    }
    
    \bib{feigin2011yangians}{article}{
          author={Feigin, Boris},
          author={Finkelberg, Michael},
          author={Negut, Andrei},
          author={Rybnikov, Leonid},
           title={Yangians and cohomology rings of laumon spaces},
            date={2011},
         journal={Selecta Mathematica},
          volume={17},
          number={3},
           pages={573\ndash 607},
    }
    
    \bib{finkelberg2000parabolic}{article}{
          author={Finkelberg, Michael},
          author={Kuznetsov, Alexander},
           title={Parabolic sheaves on surfaces and affine lie algebra},
            date={2000},
    }
    
    \bib{finkelberg2019shifted}{article}{
          author={Finkelberg, Michael},
          author={Tsymbaliuk, Alexander},
           title={Shifted quantum affine algebras: integral forms in type a},
            date={2019},
         journal={Arnold Mathematical Journal},
          volume={5},
          number={2-3},
           pages={197\ndash 283},
    }
    
    \bib{hilburn2023bfn}{article}{
          author={Hilburn, Justin},
          author={Kamnitzer, Joel},
          author={Weekes, Alex},
           title={Bfn springer theory},
            date={2023},
         journal={Communications in Mathematical Physics},
           pages={1\ndash 68},
    }
    
    \bib{laumon1988analogue}{article}{
          author={Laumon, G{\'e}rard},
           title={Un analogue global du c{\^o}ne nilpotent},
            date={1988},
    }
    
    \bib{maulik2012quantum}{article}{
          author={Maulik, Davesh},
          author={Okounkov, Andrei},
           title={Quantum groups and quantum cohomology},
            date={2012},
         journal={arXiv preprint arXiv:1211.1287},
    }
    
    \bib{nakajima2011handsaw}{article}{
          author={Nakajima, Hiraku},
           title={Handsaw quiver varieties and finite w-algebras},
            date={2011},
         journal={arXiv preprint arXiv:1107.5073},
    }
    
    \bib{neguct2013quantum}{article}{
          author={Negu{\c{t}}, Andrei},
           title={Quantum toroidal and shuffle algebras},
            date={2013},
         journal={arXiv preprint arXiv:1302.6202},
    }
    
    \bib{neguct2018affine}{article}{
          author={Negu{\c{t}}, Andrei},
           title={\space Affine laumon spaces and a conjecture of kuznetsov},
            date={2018},
         journal={\space arXiv \space preprint arXiv:1811.01011},
    }
    
    \bib{neguct2019pbw}{article}{
          author={Negu{\c{t}}, Andrei},
           title={The pbw basis of ${U_{q,p}(\ddot{\fgl}_n)}$},
            date={2019},
         journal={arXiv preprint arXiv:1905.06277},
    }
    
    \bib{okounkov2015lectures}{article}{
          author={Okounkov, Andrei},
           title={Lectures on k-theoretic computations in enumerative geometry},
            date={2015},
         journal={arXiv preprint arXiv:1512.07363},
    }
    
    \bib{okounkov2021inductive}{article}{
          author={Okounkov, Andrei},
           title={Inductive construction of stable envelopes},
            date={2021},
         journal={Letters in Mathematical Physics},
          volume={111},
           pages={1\ndash 56},
    }
    
    \bib{tsymbaliuk2010quantum}{article}{
          author={Tsymbaliuk, Aleksander},
           title={Quantum affine gelfand--tsetlin bases and quantum toroidal
      algebra via k-theory of affine laumon spaces},
            date={2010},
         journal={Selecta Mathematica},
          volume={16},
          number={2},
           pages={173\ndash 200},
    }
    
    \end{biblist}
    \end{bibdiv}
\end{document}